\documentclass{birkjour}

\usepackage{color,euscript}
\usepackage{LuminyDef}
\usepackage{epsfig,graphicx}

 \newtheorem{thm}{Theorem}[section]
 \newtheorem{alg}{Algorithm}[section]

 \theoremstyle{definition}
 
 \theoremstyle{remark}
 \newtheorem{rem}[thm]{Remark}
 
 \numberwithin{equation}{section}

\begin{document}

%-------------------------------------------------------------------------
% editorial commands: to be inserted by the editorial office
%
%\firstpage{1} \volume{228} \Copyrightyear{2004} \DOI{003-0001}
%\seriesextraline{This is the Concrete Title of this Book\br H.E. R and S.T.C. W, Eds.}
%\firstpage{1}
%\submitted{February  2014}

\title[Interpolatory Model Reduction]
 {Model Reduction by Rational Interpolation}

%----------Author 1
\author[Beattie]{Christopher Beattie}

\address{%
Department of Mathematics\\
Virginia Tech.\\
Blacksburg, VA, 24061-0123\\
USA}

\email{beattie@math.vt.edu}

\thanks{This work was supported in part by the National
Science Foundation under contract DMS-1217156.}
%----------Author 2
\author[Gugercin]{Serkan Gugercin}

\address{%
Department of Mathematics\\
Virginia Tech.\\
Blacksburg, VA, 24061-0123\\
USA}

\email{gugercin@math.vt.edu}

%%----------Author 2
%\author{A Second Author}
%\address{The address of\br
%the second author\br
%sitting somewhere\br
%in the world}
%\email{dont@know.who.knows}
%%----------classification, keywords, date
%\subjclass{Primary 99Z99; Secondary 00A00}
%
%\keywords{Class file, journal}


\keywords{Rational interpolation, model reduction, $\htwo$ approximation, 
parametric systems, generalized coprime factorization, weighted model reduction, Loewner framework 
\\ \\   \textsc{ Submitted for publication in February 2014}
}

\date{February 28, 2014}

\begin{abstract}
The last two decades have seen major developments in interpolatory methods for model reduction of large-scale linear dynamical systems.   
Advances of note include the ability to 
produce (locally) optimal reduced models at modest cost; refined methods for deriving interpolatory reduced models directly from input/output measurements; and extensions for the reduction of  parametrized systems. This chapter offers a survey of interpolatory model reduction methods starting from basic principles and ranging up through recent developments that include weighted model reduction and structure-preserving methods  based on 
 generalized coprime representations.  Our discussion is supported by an assortment of numerical examples.
\end{abstract}

%%% ----------------------------------------------------------------------
\maketitle

\vspace{-0.75cm}
\section{Introduction}
Numerous techniques exist for model reduction of large-scale dynamical systems among them, Proper Orthogonal Decomposition (see Chapter \ref{chap:Volkwein}), Balanced Truncation (see Chapter \ref{chap:Breiten}), and Interpolatory Methods, to be discussed both here and in Chapter \ref{chap:Antoulas}.   Interpolatory model reduction methods include methods referred to as Rational Krylov Methods, but should be viewed as distinct for reasons we describe later. Over the past two decades, major progress has been made in interpolation-based model reduction approaches and as a result, these methods have emerged as one of the leading choices for reducing large scale dynamical systems.

This chapter gives a survey of projection-based interpolatory methods for model reduction.  
Section \ref{sec:setting} introduces the model reduction problem setting that we consider. In Section \ref{sec:intplt},  we give general projection results for interpolatory model reduction followed by a discussion of 
$\mathcal{H}_2$ optimal model reduction by interpolation which is presented in Section \ref{sec:h2}. 
Up through Section \ref{sec:h2}, we assume  that the original system to be reduced is in a standard first-order state-space form.  Beginning in Section \ref{sec:coprime}, we discuss how  interpolatory methods can be extended with ease to much more general settings that include  systems with delays and systems with polynomial structure.   For the most part, we assume that the internal dynamics of the full order system are specified, accessible, and manifested either in a known state-space  generalized coprime representation.  This assumption will be relaxed in Section \ref{sec:tfirka} (and also in Chapter \ref{chap:Antoulas}), where a data-driven framework for interpolatory methods is introduced, useful for situations  with no direct access to internal dynamics and where only input/output measurements are available.    Finally in Section \ref{sec:pmor}, we show how to use interpolatory methods for reducing parametric dynamical systems.

The methods we discuss in this chapter have been applied with great success to very large-scale dynamical systems. Motivated by brevity, we do not present such examples here, preferring instead to illustrate the ideas with simple and approachable (albeit more academic) examples that may better reveal details of the process. For those with a hunger for more, we refer to the original papers where large-scale case studies are presented. 

\section{Model Reduction via Projection} \label{sec:setting} 
In this section, we introduce the basic concepts of projection-based model reduction. 
We also discuss the main error measures with which the approximation error will be quantified. 
\subsection{The problem setting}
We consider linear dynamical systems represented in state-space form as:
\begin{equation} \label{ltisystemintro}
  \begin{array}{rcl}
       \bE \dot{\bx}(t) & = & \bA\bx(t)+\bB\bu(t)  \\
       \by(t) & = & \bC\bx(t) + \bD \bu(t)
    \end{array}\quad \mbox{with }~\bx(0)=\mathbf{0},
\end{equation}
where $\bA,\bE\in {\mathbb R}^{n\times n}$, $\bB\in {\mathbb
R}^{n\times m}$,  $\bC\in {\mathbb R}^{p\times n}$, and $\bD\in {\mathbb R}^{p\times m}$ are 
constant matrices. In (\ref{ltisystemintro}), 
 $\bx(t) \in {\mathbb R}^n$ is the {\em internal variable}, or the {\em state variable} if $\bE$ is non-singular. 
 The length, $n$, of $\bx$ is called the {\em dimension} of the underlying dynamical system.  
 $\bu(t) \in {\mathbb R}^m$
 and $\by(t) \in {\mathbb R}^p$ are, respectively,  
the {\em inputs} and   {\em outputs} of the system.  Dynamical systems 
with $m=p=1$ will be called SISO systems ({\em single-input (and) single-output})
  while all other cases will be grouped together and referred to as 
 MIMO systems ({\em multi-input (or) multi-output}). 

For cases where the  dimension $n$ is large, e.g., $ n \geq 10^5,10^6$, the 
simulation and control of the  system can lead to a huge computational burden; 
especially when the system must be resimulated over and over again, 
say, using different input selections, $\bu(t)$. 
Our goal is to replace (\ref{ltisystemintro}) with a simpler 
reduced model having the form 
\begin{equation} \label{redsysintro}
 \begin{array}{rcl}
   \bE_r\dot{\bx}_{r}(t) & = & \bA_{r}\bx_{r}(t)+\bB_{r}\bu(t)  \\
   \by_{r}(t) & = & \bC_{r}\bx_{r}(t)  + \bD_r \bu(t)
\end{array}\quad \mbox{with }~\bx_{r}(0)=\mathbf{0},
\end{equation}
	where  $\bA_r,\bE_r \in{\mathbb R}^{r \times r}$, $\bB_r \in {\mathbb
R}^{r \times m}$, $\bC_r \in{\mathbb R}^{p \times r}$ and $\bD_r\in {\mathbb R}^{p\times m}$
with $r \ll n$, and such that over a wide range of system inputs, the corresponding outputs
 of the reduced system, $\by_r(t)$, will be good approximations to the corresponding true outputs,
$\by(t)$, with respect to an appropriate, physically relevant norm.
\subsection{Transfer Function and Error Measures}
For the linear dynamical systems considered here, the frequency domain representation  is a powerful tool to quantify the model reduction error. Indeed, in this setting,  error analysis in the frequency domain relates directly  to  error analysis in the time domain.

Let $\widehat{\by}(s)$, $\widehat{\by}_r(s)$, and $\widehat{\bu}(s)$ denote the Laplace transforms of $\by(t)$, $\by_r(t)$ and $\bu(t)$, respectively.  Taking the Laplace transforms of (\ref{ltisystemintro}) and (\ref{redsysintro}) yields 
\begin{align}
 \widehat{\by}(s)&=\left(\bC \left(s \bE -\bA\right)^{-1}\bB+\bD\right)\,\widehat{\bu}(s),
   \label{LaplTrLTI}\\
  \widehat{\by}_r(s)&=\left(\bC_{r} \left(s \bE_r-\bA_{r}\right)^{-1}\bB_{r}+\bD_r\right)\,\widehat{\bu}(s).
   \label{LaplTrLTIred}
\end{align}
The mapping from $\widehat{\bu}(s)$ to $\widehat{\by}(s)$ is called the {\em transfer function}. 
Likewise, the mapping from $\widehat{\bu}(s)$ to $\widehat{\by}_r(s)$  is the {\em transfer function of the reduced model}. We denote them by $\bH(s)$ and $\bH_r(s)$, respectively:
\begin{align} 
\bH(s)=& ~\bC \left(s \bE -\bA\right)^{-1}\bB+ \bD, \mbox{ and } \label{TransFnc} \\
\bH_r(s)=&~\bC_r \left(s \bE_r -\bA_r\right)^{-1}\bB_r+ \bD_r.  \label{redTransFnc}
\end{align}
$\bH(s)$ is a $p \times m$ matrix-valued degree-$n$ rational function in $s$, and
$\bH_r(s)$ is a $p \times m$ matrix-valued degree-$r$ rational function in $s$. 
Thus, the model reduction problem could be viewed as a rational approximation problem in the complex domain. This perspective of model reduction is emphasized in Section \ref{sec:intplt}.

\subsection{Petrov-Galerkin Projections} Most model reduction methods can be formulated with the aid of either Petrov-Galerkin or Galerkin projections.  Even though the original internal variable, $\bx(t)$, evolves in a (large) 
$n$-dimensional space, it is often the case that it hews rather closely to some (typically unknown) $r$ dimensional subspace. Let $\bV \in \IR^{n \times r}$ be a basis for this subspace, which is as yet undetermined. The original state may be approximated as  $\bx(t) \approx \bV \bx_r(t)$ for some $\bx_r(t) \in \IR^r$ and this expression may be used to represent the reduced model dynamics.  Plug the approximation, $\bx(t) \approx \bV \bx_r(t)$,  into (\ref{ltisystemintro}) to obtain a residual 
\begin{equation} \label{eq:res}
\bR(\bx_r(t)) =  \bE \bV\dot{\bx}_r (t)-\bA\bV\bx_r (t)-\bB\,\bu(t)
\end{equation}
and the approximate output 
\begin{equation}\label{eq:yr}
\by_r(t) = \bC \bV \bx_r(t)  + \bD \bu(t).
\end{equation}
The reduced state trajectory, $\bx_r(t)$, is determined by enforcing a Petrov-Galerkin orthogonality condition on the residual $\bR(t)$ in (\ref{eq:res}): we pick another $r$ dimensional subspace with a basis $\bW \in \IR^{n \times r}$ and impose the Petrov-Galerkin condition:
$$
\bW^T \bR(t) = \bW^T \left( \bE \bV\dot{\bx}_r (t)-\bA\bV\bx_r (t)-\bB\,\bu(t) \right) = \mathbf{0}.
$$
This leads to a reduced model as in (\ref{redsysintro}),
$$
   \bE_r\dot{\bx}_{r}(t)  =  \bA_{r}\bx_{r}(t)+\bB_{r}\bu(t),
   \quad   \by_{r}(t)  =  \bC_{r}\bx_{r}(t)  + \bD_r \bu(t)
$$
 with reduced model quantities defined as 
 \small{
 \begin{align} \label{red_projection}
\bE_r=\bW^{T}\bE\bV,~~~
\bA_r = \bW^{T} \bA \bV,~~~ \bB_r= \bW^{T}\bB,~~~\bC_r = \bC\bV,\mbox{ and }
\bD_r = \bD.
\end{align}
}
A critical observation to make here is that the reduced model does not depend on the specific basis selection 
made for $\bV$ and $\bW$, only on the subspaces themselves. To see this, let 
$\widetilde{\bV} = \bV \bT_1$ and $\widetilde{\bW} = \bW \bT_2$ where 
$\bT_1 \in \IR^{r \times r}$ and $\bT_2 \in \IR^{r \times r}$ are nonsingular matrices corresponding to a change-of-basis. This leads to a change in reduced model quantities 
$$
\widetilde{\bE}_r=\bT_2^T\bE_r\bT_1,~~~
\widetilde{\bA}_r= \bT_2^T\bA_r\bT_1,~~~ \widetilde{\bB}_r= \bT_2^T\bB_r,~~~\widetilde{\bC}_r=  \bC_r\bT_1, \mbox{ and }
\widetilde{\bD}_r = \bD_r,
$$
where $\bE_r$, $\bA_r$, $\bB_r$, $\bC_r$, and $\bD_r$ are as defined in (\ref{red_projection}). 
A straightforward  comparison of the transfer functions reveals 
\begin{align*}
\widetilde{\bH}_r(s) &= \widetilde{\bC}_r \left(s \widetilde{\bE}_r -\widetilde{\bA}_r \right)^{-1}\widetilde{\bB}_r+ \widetilde{\bD}_r\\
&=  \bC_r\bT_1 \left(s \bT_2^T\bE_r \bT_1-\bT_2^T\bA_r\bT_1\right)^{-1}\bT_2^T\bB_r+ \bD_r \\
& = \bC_r \left(s \bE_r -\bA_r\right)^{-1}\bB_r+ \bD_r  = \bH_r(s).
\end{align*}
\subsection{Error Measures}
As in any approximation problem, error measures are necessary to quantify the approximation error appropriately. For  linear dynamical systems, error analysis is best presented in the frequency domain, yet directly relates to time-domain error,  $\by(t) - \by_r(t)$. 
It follows from (\ref{LaplTrLTI}) and (\ref{LaplTrLTIred}) that
$$
\widehat{\by}(s)-\widehat{\by}_r(s)=\left[\bH(s) - \bH_r(s)\right]\, \widehat{\bu}(s).
$$
Thus, the closeness of $\widehat{\by}_r(s)$ to $\widehat{\by}(s)$ is directly related to 
the closeness of $\bH_r(s)$ to $\bH(s)$. The $\htwo$ and $\hinf$ norms are the most common measures of closeness for transfer functions:
\subsubsection{The $\hinf$ Norm} 
Let $\bH(s)$ be the transfer function of a stable dynamical system. The $\hinf$ norm is defined 
as 
\begin{equation} \label{eq:hinf}
\left \| \bH \right\|_\hinf = \sup_{\omega \in \IR} \left\| \bH(\imath \omega) \right\|_2 ,
\end{equation}
where $\left\| \mathbf{M} \right\|_2$ denotes the spectral (Euclidean-induced) norm of the complex matrix $\bM$.
When $\bE$ is nonsingular, all eigenvalues of the matrix pencil 
$\lambda \bE - \bA$ must lie in the left-half plane. When $\bE$ is singular, we assume additionally that
$0$ is not a defective eigenvalue of $\bE$.  This guarantees that $\bH(s)$ 
remains bounded as $s \to \infty$.

The importance of the $\hinf$-norm stems from it being the ($L_2$-$L_2$)-induced operator norm of an underlying convolution operator mapping the system inputs, $\bu$,
 to system outputs, $\by$:
 \small{$\displaystyle 
\left \| \bH \right\|_\hinf  = \sup_{\bu\in L_2}\frac{ \left\| \by \right\|_{L_2}}{ \left\| \bu \right\|_{L_2}}.
$}
where $\left\| \bz \right\|_{L_2} =\sqrt{ \int_0^\infty\left\|\bz(t)\right\|_2^2\,dt } $,
With respect to model reduction error, one directly obtains 
$$
\left\| \by - \by_r\right\|_{L_2} \leq \left \| \bH - \bH_r\right\|_\hinf  \left\| \bu\right\|_{L_2}.
$$
If one wishes to produce reduced models that generate outputs, $\by_r(t)$, 
that are always close (with respect to the $L_2$ norm) to the corresponding true outputs, $\by(t)$, uniformly 
so over all $L_2$-bounded inputs, $\bu(t)$, 
then one should apply a model reduction technique that produces small $\hinf$ error. 

\subsubsection{The $\htwo$ Norm} \label{subsubsec:h2norm} 
Let $\bH(s)$ be the transfer function of a stable dynamical system. Then, the $\htwo$ norm is defined  as
\begin{eqnarray}
\left\| \bH \right\|_{\htwo} &:=&
\left(\frac{1}{2\pi}\int_{-\infty}^{\infty}
  \left\| \bH(\imath \omega )  \right\|_{\rm F}^2\right)^{1/2},
\end{eqnarray}
where  $\left\| \mathbf{M} \right\|_{\rm F}$ denotes the Frobenius norm of a complex matrix $\bM$.

When $\bE$ is nonsingular, all eigenvalues of the matrix pencil 
$\lambda \bE - \bA$ must lie in the left-half plane. When $\bE$ is singular, we assume additionally that
$0$ is not a defective eigenvalue of $\bE$.

When $\bE$ is nonsingular, we require that all the eigenvalues of the matrix pencil 
$\lambda \bE - \bA$ lie in the left-half plane and $\bD = \mathbf{0}$ for the $\htwo$ norm to be finite. 
When $\bE$ is singular, we assume in addition that
$0$ is not a defective eigenvalue of $\bE$ and that ${\displaystyle \lim_{s \to \infty} \bH(s) = \mathbf{0}}$.
The $\htwo$ norm bears a direct relationship to the time domain norm of $\by(t)$: 
$$
\left \| \by \right \|_{L_\infty} = \sup_{t>0}\, \|\by(t)\|_{\infty}  \leq \left \| \bH \right \|_{\htwo}
\left\| \bu \right\|_{L_2},
$$ 
and this bound is best possible for MISO systems ($p=1$), SIMO systems ($m=1$), and SISO  systems ($m=p=1$), reflecting the fact that 
 the $\htwo$ norm is simply the ($L_2$-$L_{\infty}$)-induced norm of the underlying convolution operator in these cases.   With respect to model reduction error, 
  we have in general, 
$$
\left \| \by - \by_r \right \|_{L_\infty}   \leq \left \| \bH - \bH_r \right \|_{\htwo}
\left\| \bu \right\|_{L_2}. 
$$ 
So, if one wishes to produce reduced models that generate outputs, $\by_r(t)$, 
that are uniformly and instantaneously close to the corresponding true outputs, $\by(t)$, 
 uniformly  so over all $L_2$-bounded inputs, $\bu(t)$, 
then one should apply a model reduction technique that produces small $\htwo$ error. 
Note that the $\htwo$ error may be made small even for original systems $\bH$ with ${\displaystyle \lim_{s \to \infty} \bH(s) \neq \mathbf{0}}$ 
with an appropriate choice of $\bD_r$ and $\bE_r$; $\left \| \bH - \bH_r \right \|_{\htwo}$ may be small even if $\left \| \bH  \right \|_{\htwo}$ is unboundedly large.   This is discussed in more detail in Section \ref{DrNotD}. 
\section{Model Reduction by Interpolation} \label{sec:intplt}
In this section, we present the fundamental projection theorems used in interpolatory model reduction and illustrate their use with some simple examples.
\subsection{Tangential Interpolation Problem}
One easily observes that the model reduction problem for linear, time invariant dynamical systems can be formulated as a rational approximation problem:  Given a degree-$n$ rational function, $\bH(s)$ (the full model), find a degree-$r$ rational function, $\bH_r(s)$ (the reduced model), that approximates $\bH(s)$ accurately with respect to either the $\hinf$ or $\htwo$ norm.  Interpolation is a commonly applied tool for function approximation; typically, effective polynomial interpolants are easy to calculate.  Here, we will develop straightforward methods for obtaining rational interpolants, however, the precise notion of ``interpolation" that will be used must be clarified. 
Since $\bH(s)$ is a $p \times m$ matrix valued rational function, the immediate extension of point-wise interpolation to matrix-valued functions suggests that we attempt to enforce  conditions such as $\bH(s_0) = \bH_r(s_0)$ at each  interpolation point $s_0 \in \IC$.  But viewed element-wise, this would require in effect, $p \times m$ interpolation conditions at every interpolation point. For systems with even modestly large number of input and outputs dimensions $m$ and $p$, this will lead to a large number of interpolation conditions requiring as a result, quite a large reduced order $r$. Thus, for MIMO systems, instead of this notion of  ``full matrix interpolation", we require only that the interpolating matrix function match the original only along certain directions, ``tangential interpolation."
We will show later that this relaxed notion of interpolation is adequate to characterize
 necessary conditions for optimal approximation in the $\htwo$ norm.

Tangential interpolation involves choosing interpolation directions in addition to interpolation points. 
We separate the interpolation points and directions into two categories: ``left" and ``right." 
We say that $\bH_r(s)$ is a right-tangential interpolant to $\bH(s)$ at $s=\sigma_i$ along the right tangent direction  $\rdir_i \in \IC^m$ if 
$$
\bH(\sigma_i) \rdir_i = \bH_r(\sigma_i) \rdir_i.
$$
Similarly, we  say that $\bH_r(s)$ is a left-tangential interpolant to $\bH(s)$ at $s=\mu_i$ along the left tangent direction $\ldir_i \in \IC^p$ if 
$$
\ldir_i^T \bH(\mu_i)= \ldir_i^T\bH_r(\mu_i).
$$
 Our model reduction task can now be formulated as tangential interpolation as follows:
Given a set of $r$ right interpolation points $\{\sigma_i\}_{i=1}^r \in \IC$, $r$ left interpolation points 
$\{\mu_i\}_{i=1}^r$, $r$ right-tangential directions  $\{\rdir_i\}_{i=1}^r \in \IC^m$,  and 
$r$ left-tangential directions  $\{\ldir_i\}_{i=1}^r \in \IC^p$, find a degree-$r$ reduced transfer function $\bH_r(s)$
so that
\begin{eqnarray}  \label{eq:l}
\begin{array}{rcl}
\bH(\sigma_i) \rdir_i &= &\bH_r(\sigma_i) \rdir_i, \\[.1in]
\ldir_i^T \bH(\mu_i) & = & \ldir_i^T\bH_r(\mu_i),
\end{array}
\quad\mbox{for}~~~ i=1,\ldots,r.
\end{eqnarray}
We  say that $\bH_r(s)$ is a bitangential Hermite interpolant to $\bH(s)$ at $s=\sigma_i$ along
 the right tangent direction $\rdir_i \in \IC^m$ and 
 the left tangent direction $\ldir_i \in \IC^p$, if 
$$
\ldir_i^T \bH'(\sigma_i) \rdir_i= \ldir_i^T\bH'_r(\sigma_i) \rdir_i.
$$
where $'$ denotes differentiation with respect to $s$. Therefore, 
in addition to (\ref{eq:l}),
we may also require $\bH_r(s)$ to satisfy
\begin{equation} \label{eq:h}
\ldir_i^T \bH'(\sigma_i) \rdir_i= \ldir_i^T\bH'_r(\sigma_i) \rdir_i, \quad\mbox{for}~~~ i=1,\ldots,r.
\end{equation}
In Section \ref{sec:h2}, we will show how to choose interpolation points and tangent directions to produce optimal approximation with respect to the $\htwo$ norm.
 \subsection{Petrov-Galerkin Projections for Tangential Interpolation}
 Our goal here is to appropriately pick the model reduction bases $\bV$ and $\bW$ so that the reduced model obtained by a Petrov-Galerkin projection as in (\ref{red_projection}) satisfies the tangential interpolation conditions (\ref{eq:l}) and (\ref{eq:h}). We first present the projection theorem 
 followed by a historical perspective:
  \begin{thm}
\label{thm:interpolation}
Given the transfer function, $\bH(s) = \bC(s\bE-\bA)^{-1}\bB + \bD$,  
let $\bH_r(s)$ denote a reduced transfer function obtained by projection as in $($\ref{red_projection}$)$ using the model reduction bases $\bV$ and $\bW$. 
For interpolation points $\sigma,\,\mu \in \IC$, suppose that $\sigma\,\bE\, -\,\bA$ and $\mu\,\bE\, -\,\bA$ be invertible.
 Let $\rdir\in \IC^{m}$ and  $\ldir\in \IC^{\ell}$  be designated (nontrivial)  tangent directions. Then,  
      \begin{itemize}
        \item[(a)]~if 
	\begin{equation}  \label{eq:rtc}
	\left(\sigma\,\bE-\bA\right)^{-1}\bB\rdir\in {\mbox{\normalfont\textsf{Ran}}}(\bV),
\end{equation}
	  then 
	 	\begin{equation}  \label{eq:rti}
	 \bH(\sigma)\rdir=\bH_{r}(\sigma)\rdir;
	\end{equation}
        \item[(b)]~if 
	\begin{equation}  \label{eq:ltc}
	 \left(\ldir^T\bC\left(\mu\,\bE -\bA\right)^{-1}\right)^{T}\in\mbox{\normalfont\textsf{Ran}}(\bW),
	 \end{equation}
	  then 
	  	\begin{equation}  \label{eq:lti}
	 \ldir^T\bH(\mu)=\ldir^T\bH_{r}(\mu);
	 	 \end{equation}
     \item[(c)]~if both $($\ref{eq:rtc}$)$ and $($\ref{eq:ltc}$)$ hold, and $\sigma = \mu$,    then 
	 	\begin{equation}  \label{eq:bti}
\ldir^T\bH'(\sigma) \rdir =\ldir^T\bH_{r}'(\sigma) \rdir
\end{equation}
 as well. 
    \end{itemize}
\end{thm}

  Theorem \ref{thm:interpolation} illustrates that imposing either a left or right tangential interpolation condition requires adding only one vector either to the left or right model reduction bases.  
  For the case of repeated left and right interpolation points, 
 the bitangential Hermite condition is satisfied for free, in the sense that no additional vectors must be included in the model reduction bases. Notice that the values that are interpolated are never explicitly computed; this is a significant advantage of the Petrov-Galerkin projection framework as used in interpolatory model reduction.

A projection framework for interpolatory model reduction was  
introduced by Skelton
{\it et. al.} in \cite{devillemagne1987mru, yousouff1984cer, yousuff1985lsa}. 
This approach was put into a robust numerical framework
by Grimme \cite{grimme1997krylov}, who employed the rational Krylov 
subspace method of Ruhe \cite{ruhe1984rational}.  
The tangential interpolation framework of Theorem \ref{thm:interpolation}
was developed by Gallivan {\it et al.}  \cite{gallivan2005mrm}. 
For SISO systems, the model reduction bases, $\bV$ and $\bW$, produced from Theorem \ref{thm:interpolation} become rational Krylov subspaces; and so interpolatory model reduction is sometimes referred to as {\em rational Krylov methods}. However, the connection to rational Krylov subspaces is lost for general MIMO systems, unless all  tangent directions are the same.  So, we prefer the simpler descriptive label, {\em interpolatory methods}. Indeed, for the more general systems that we consider in Section \ref{sec:coprime}, the direct extensions we develop for interpolatory methods have no connection to rational Krylov subspaces even in the SISO case.  Another term that has been in use when referring to interpolatory model reduction methods is {\em Moment Matching Methods}. The $k{\rm th}$ moment of $\bH(s)$ around $\sigma$ is the $k{\rm th}$ derivative of $\bH(s)$ evaluated at $s =\sigma$. 
In the SISO case, the reduced transfer function obtained via rational interpolation will match those moments - in effect, generalized Hermite interpolation.  The notion of moment matching for MIMO systems with respect to tangent directions is not so clearly stated, however.  
See \cite{Ant2010imr,FelF95, GalGV94,Bai02,Fre03,antoulas2005approximation,BenHtM11} and the references therein for other related work on model reduction by interpolation.

\subsubsection*{A simple example} Consider the following linear dynamical system with $n=3$,  $m=p=2$, $\bE = \bI_3$, $\bD = \mathbf{0}$, 
$$
\bA = \left[ \begin{array}{rrr}
-6 & -11 & -6 \\ 1& 0& 0 \\ 0 & 1 &0 \end{array} \right],\quad
\bB = \left[ \begin{array}{rr}
-1 &1 \\ 0 & 1 \\ 1 & 0\end{array} \right],\quad\mbox{and}\quad
\bC =  \left[\begin{array}{rrr}
1 & 0 & 1 \\ 1& -1 & 0 \end{array} \right].
$$
The transfer function of this dynamical system can be computed as:
$$
\bH(s) =
\frac{1}{s^3 + 6 s^2 + 11 s + 6}
\left[\begin{array}{cc}
10 &  s^2 - 10 s + 1 \\
-s^2 - 5 s + 6 &-18 s - 6 
\end{array} \right].
$$
Let $\sigma_1 = \mu_1 = 0$ be the left and right interpolation points together with 
tangent directions 
$$
\rdir_1 = \left[\begin{array}{c}
1 \\ 2 \end{array} \right]  
\quad\mbox{and}\quad
\ldir_1 = \left[\begin{array}{c}
3 \\ 1 \end{array} \right].
$$
Using Theorem \ref{thm:interpolation}, we compute the interpolatory model reduction bases:
\begin{align*}
\bV &= (\sigma_1\bE -\bA)^{-1}\bB\rdir_1 =   \left[\begin{array}{r}
-2 \\ -1 \\ 4
\end{array} \right] \mbox{ and } \\
\bW &= (\sigma_1\bE -\bA)^{-T}\bC^T\ldir_1 =  \left[\begin{array}{r}
0.5 \\ -1 \\ 6.5 
\end{array} \right].
\end{align*}
Then, the Petrov-Galerkin projection  in (\ref{red_projection}) leads to the reduced model quantities
$$
\bE_r = 26, ~\bA_r = -5,~ \bB_r =  \left[\begin{array}{rr}
6 & -0.5
\end{array} \right],~
\bC_r = 
 \left[\begin{array}{r}
 2 \\ -1
\end{array} \right],~\bD_r = \mathbf{0},
$$
and consequently to the reduced model transfer function
$$
\bH_r(s) = \frac{1}{26 s + 5} 
 \left[\begin{array}{rr}
 12 & -1 \\
 -6 & 0.5
\end{array} \right].
$$
%(only four significant digits are printed). 
Now that we have $\bH_r(s)$, we can check the interpolation conditions explicitly:
\begin{align*}
\bH(\sigma_1)\rdir_1 &= \bH_r(\sigma_1)\rdir_1=  \left[\begin{array}{r}
2 \\ -1
\end{array} \right],   \quad \checkmark  \\
\ldir_1^T \bH(\sigma_1)& = \ldir_1^T \bH_r(\sigma_1) =  \left[\begin{array}{rr}
6 & -0.5
\end{array} \right],\quad \checkmark  \\
\ldir_1^T \bH'(\sigma_1) \rdir_1 & = \ldir_1^T \bH_r'(\sigma_1) \rdir_1 = -26. \quad \checkmark 
\end{align*}
Notice that $\bH_r(s)$ does not fully interpolate $\bH(s)$ at $s=\sigma_1=0$:
$$
\bH(\sigma_1) = \left[\begin{array}{cc}
    5/3 &   1/6 \\
    1   &  -1
\end{array} \right] \neq 
\left[\begin{array}{rr}
    2.4 &   -0.2 \\
    -1.2   &  0.1
\end{array} \right] = \bH_r(\sigma_1).
$$
 To enforce full matrix interpolation, we need to modify the construction of the model reduction bases and remove the tangential vectors. Denote the new bases for full matrix interpolation by 
 $\bV_{\! \!m}$ and $\bW_{\! \! m}$:
\begin{align*}
\bV_{\! \!m} &= (\sigma_1\bE -\bA)^{-1}\bB =   \left[\begin{array}{rr}
0 & -1 \\ -1 & 0  \\ 5/3 & 7/6 
\end{array} \right] \\
\bW_{\! \! m} &= (\sigma_1\bE -\bA)^{-T}\bC^T =  \left[\begin{array}{rr}
1/6 & 0 \\ 0 & -1 \\ 
11/6 & 1
\end{array} \right].
\end{align*}
Note that even with a single interpolation point, in the case of full matrix interpolation, the reduction spaces has dimension $2$ leading a degree-$2$ reduced model with transfer function
%$$
%\bH_{rm}(s) = \frac{1}{s^2 + 2.082 s + 1.082}
% \left[\begin{array}{cc}
% 0.1639 s + 1.803   &  -1.787 s + 0.1803 \\
% -0.8033 s + 1.082 & -3.344 s - 1.082
%\end{array} \right].
%$$
$$
\bH_{rm}(s) =\left[\begin{array}{cc} 10 & 1\\  6 &  -6 \end{array} \right]
\left(s\left[\begin{array}{cc} 110 & 71\\  96 &  42 \end{array} \right] -
\left[\begin{array}{cc} -60 & -6 \\  -36 &  36 \end{array} \right]\right)^{-1}
\left[\begin{array}{cc} 10 & 1\\  6 &  -6 \end{array} \right]
$$
 This new reduced model fully interpolates $\bH(s)$ and $\bH'(s)$ at $s=\sigma_1 =0$:
 \begin{align*}
\bH(\sigma_1) &= \bH_{rm}(\sigma_1) = \left[\begin{array}{cc}
    5/3 &   1/6 \\
    1   &  -1
\end{array} \right] \\
 \bH'(\sigma_1) &= \bH'_{rm}(\sigma_1) = \left[\begin{array}{cc}
     -55/18 &   -71/36\\
 -8/3   &  -7/6
\end{array} \right].  \qquad\qquad 
\end{align*}    
 This simple example illustrates the fundamental difference between tangential and full interpolation. In the case of tangential interpolation, each interpolation condition only adds one degree of freedom to the reduced dimension; however in the case of full interpolation, a single (full matrix) interpolation condition will generically add $m$ or $p$ dimensions to the reduced model.  We will see in \S \ref{sec:h2} that optimality requires only tangential interpolation. \quad
 $\Box$
 
%$ 
% \left[\begin{array}{r}
%\end{array} \right].
%$
The computations in the previous example can be extended easily to the case of 
$r$ interpolation points: Given $\bH(s) = \bC(s\bE - \bA)^{-1}\bB + \bD$,
$r$  right interpolation points $\{\sigma_i\}_{i=1}^r$ and 
right  directions $\{ \rdir_k\}_{k=1}^r \in \IC^m$, 
 $r$  left interpolation points $\{\mu_j\}_{j=1}^r$ and
left  directions  $\{ \ldir_k\}_{k=1}^r \in \IC^p$, construct 
\begin{align}
\bV =&
\left[(\sigma_1\bE-\bA)^{-1}\bB\rdir_1,~\cdots,~
(\sigma_r\bE-\bA)^{-1}\bB\rdir_r\right] \ \mbox{and} \label{eqn:Vr} 
\end{align}
\begin{align}
\bW =&
\left[(\mu_1\bE-\bA)^{-T}\bC^T\ldir_1,~\cdots,~
(\mu_r\bE-\bA)^{-T}\bC^T\ldir_r\right] \ \label{eqn:Wr}.
\end{align}
Then,  $\bH_r(s)= \bC_r (s\bE_r - \bA_r)^{-1}\bB_r$ 
constructed by a Petrov-Galerkin projection as in (\ref{red_projection})
satisfies the Lagrange tangential interpolation conditions (\ref{eq:l}) and 
 the bitangential Hermite interpolation conditions (\ref{eq:h}), if in addition $\sigma_i = \mu_i$
 (provided 
that $\sigma_i\bE_r-\bA_r$  and $\mu_i\bE_r-\bA_r$ are nonsingular for each 
$i=1,\cdots,\, r$).  

 Theorem \ref{thm:interpolation} can be  extended readily to include higher-order Hermite interpolation: 

 \begin{thm}
\label{thm:interpolation_highorder}
Given a full order model with transfer function $$\bH(s) = \bC(s\bE-\bA)^{-1}\bB + \bD,$$
 let $\bH_r(s)$ denote a reduced transfer function obtained by projection as in $($\ref{red_projection}$)$,
  using model reduction bases, $\bV$ and $\bW$. 
Let $\bH^{(k)}(\sigma)$ denote the $k{\rm th}$ derivative of $\bH(s)$ with respect to $s$ evaluated at $s=\sigma$.
For interpolation points $\sigma,\,\mu \in \IC$, suppose $\sigma\,\bE\, -\,\bA$ and $\mu\,\bE\, -\,\bA$ are invertible,
and that $\rdir\in \IC^{m}$ and  $\ldir\in \IC^{p}$  are  given (nontrivial)  tangent directions. Then,  
      \begin{itemize}
        \item[(a)]~if 
	\begin{equation}  \label{eq:hortc}
	\left(\left(\sigma\,\bE-\bA\right)^{-1}\bE \right)^{j-1}\left(\sigma\,\bE-\bA\right)^{-1} \bB\rdir\in  \mbox{\normalfont\textsf{Ran}}(\bV),
	\mbox{ for } j=1,.,N  
\end{equation}
	  then 
	 	\begin{equation}  \label{eq:horti}
	\bH^{(k)}(\sigma)\rdir=\bH_r^{(k)}(\sigma)\rdir~~for~~k = 0,1,\ldots,N-1;
	\end{equation}
        \item[(b)]~if 
	\begin{equation}  \label{eq:holtc}
	\left(\left(\mu\,\bE-\bA\right)^{-T}\bE^T \right)^{j-1}\left(\mu\,\bE-\bA\right)^{-T} \bC^T\ldir\in  \mbox{\normalfont\textsf{Ran}}(\bW)
	 \mbox{ for }j=1,.,M,
	 \end{equation}
	  then 
	  	\begin{equation}  \label{eq:holti}
	\ldir^T\bH^{(k)}(\mu)=\ldir^T\bH_r^{(k)}(\mu)~~ for~~k = 0,1,\ldots,M-1;
	 	 \end{equation}
		 
		 \vspace{-.1ex}
     \item[(c)]~if $\sigma = \mu$ and both $($\ref{eq:hortc}$)$ and $($\ref{eq:holtc}$)$ hold,    
      then 
	 	\begin{equation}  \label{eq:hobti}
\ldir^T\bH^{(k)}(\sigma) \rdir =\ldir^T\bH^{(k)}_{r}(\sigma) \rdir
~~ for~~k = 0,\ldots,M+N-1
\end{equation}
 as well. 
    \end{itemize}
\end{thm}

The main cost in interpolatory model reduction originates from the need to solve large-scale (typically sparse) shifted linear systems.  There is no need to solve large-scale Lyapunov or Riccati equations,  giving interpolatory methods a computational  advantage over competing methods.  The discussion here assumes that these linear systems are solved by direct methods (e.g., Gauss elimination). However, for systems with millions of degrees of freedom, one would prefer to incorporate iterative solution strategies to construct the model reduction bases, $\bV$ and $\bW$.  We refer to \cite{beattie2012inexact,wyatt2012issues,beattie2006inexact} for detailed analyses of the effects of iterative solves on interpolatory model reduction and to 
\cite{beattie2012inexact,ahuja2012recycling,ahmad2012preconditioned,benner2011recycling,MPIMD13-21} for development of effective iterative solves in the context of interpolatory model reduction.

\subsubsection{Rational interpolants with  $\bD_r \neq \bD$} \label{DrNotD}
So far, we have assumed that $\bD_r = \bD$. This is the logical choice if one is interested in minimizing the $\htwo$ norm of the error system. For the case of ordinary differential equations where $\bE$ is nonsingular, choosing $\bD_r \neq \bD$ will lead to unbounded $\htwo$ error norm. However,  if, instead,  one is interested in  $\mathcal{H}_{\infty}$ error, then  flexibility in choosing  $\bD_r $ will be  necessary as the optimal $\hinf$ approximation will have $\bD_r\neq \bD$ (see, e.g., \cite{beattie2008iba,flagg2013interpolatory}).

Another case that may require choosing  $\bD_r\neq \bD$ is the case of an index 1 system of differential algebraic equations (DAEs). In our setting,  this means that the $\bE$ matrix in   (\ref{ltisystemintro}) has a non-defective eigenvalue at $0$.  (Interpolatory projection methods for DAEs is considered in detail in Section \ref{sec:dae} below.) 
In this case, $\displaystyle \lim_{s\rightarrow\infty} \bH(s) \neq  \bD$, 
so for $\bH_r(s)$ to match $\bH(s)$ asymptotically well at high frequencies,  we require 
$$
\bD_r=\lim_{s\rightarrow\infty} \left(\bH(s)-\bC_r(s\bE_r-\bA_r)^{-1}\bB_r\right).
$$
Since $\bE_r$ will be generically nonsingular (assuming $r<\mbox{rank}(\bE)$),
setting   $\bD_r =  \lim_{s\rightarrow\infty} \bH(s)$ will guarantee  $\lim_{s\rightarrow\infty} \bH(s) = 
\lim_{s\rightarrow\infty} \bH_r(s)$. 

The next theorem shows how one may construct reduced-models with  $\bD_r\neq \bD$ without losing interpolation properties.  Without loss of generality, we  assume  $\bD=\mathbf{0}$, i.e.
\begin{equation} \label{HwoD}
     \bH(s)=\bC(s\bE-\bA)^{-1}\bB.
 \end{equation}
 In the general case with $\bD \neq \mathbf{0}$, one simply need to replace $\bD_r$ with $\bD_r-\bD$.
The result below was first given in
\cite{mayo2007fsg} and later generalized in \cite{beattie2008ipm}.
 \begin{thm} \label{GenIntrp2}
 Given are a full-order model with transfer function $\bH(s)$ as in (\ref{HwoD}),
 $r$ distinct left-interpolation points  $  \{ \mu_i\}_{i=1}^r$ together with 
 $r$ left tangent directions $ \{ \ldir_i\}_{i=1}^r\subset\IC^{p} $, and,
 $r$ distinct right-interpolation points 
         $ \{ \sigma_j \}_{j=1}^r,$ together with 
 $r$ right tangent directions  $ \{ \rdir_j \}_{j=1}^r\subset\IC^{m}$. 
       Let the model reduction bases $\bV_r \in \IC^{n \times r}$ 
    and  $\bW_r\in \IC^{n \times r}$  be constructed as in  (\ref{eqn:Vr}) and (\ref{eqn:Wr}), respectively.
          Define $\widetilde{\rdir}$ and $\widetilde{\ldir}$
    as        $$
    \widetilde{\rdir} =[\rdir_1,\,
     \rdir_2,\,...,\,
     \rdir_r\,]
     \quad\mbox{and}\quad \widetilde{\ldir}^T=
     \left[
      \ldir_1,\,
      \ldir_2,\,
        \ldots,\,
       \ldir_r \right]^T
       $$
For any $\bD_{r}\in  \IC^{p \times m}$, define
\begin{align}   
\bE_r &= \bW_r^T \bE \bV_r, \qquad \bA_r = \bW_r^T 
\bA \bV_r + \widetilde{\ldir}^T \bD_r \widetilde{\rdir}, \nonumber \\
\bB_r  &= \bW_r^{T} \bB- \widetilde{\ldir}^T \bD_{r}, \qquad {\rm and}~~
\bC_r  = \bC \bV_r - \bD_r \widetilde{\rdir}
\label{ROM_DnotDr}
\end{align}
Then the reduced-order model
$
    \bH_r(s)=\bC_r(s \bE_r - \bA_r)^{-1}\bB_r+\bD_{r}
$
satisfies       
$$
      \bH(\sigma_i) \rdir_i=\bH_{r}(\sigma_i) \rdir_i  \quad \mbox{and}\quad 
          \ldir_i^T\bH(\mu_i) =\ldir_i^T\bH_{r}(\mu_i)
                   \quad\mbox{ for $i=1,\,...,\,r$. }
$$ 
\end{thm}
Theorem \ref{GenIntrp2} shows how to construct a rational tangential interpolant with an {\it arbitrary} $\bD_r$ term. This $\bD_r$  can be chosen to satisfy specific design goals. 

\subsection{Interpolatory Projections for Differential Algebraic Systems} \label{sec:dae}
The interpolation conditions in Theorems \ref{thm:interpolation} and \ref{thm:interpolation_highorder} are valid regardless of whether or not the matrix $\bE$ is singular, as long as $s \bE - \bA$ and $s \bE_r - \bA_r$ are invertible matrices for $s=\sigma,\mu$.  When $\bE$ is nonsingular, the underlying model is a system of ordinary differential equations (ODE);  when $\bE$ is singular, the underlying model is a system of  differential algebraic equations (DAE)).  Thus, from a pure interpolation perspective,  the distinction does not make a difference. However, from
the perspective of error measures, there is a crucial difference. 

A crucial difference between a DAE system and an ODE system is that the transfer function of a DAE system could contain a polynomial part that may grow unboundedly as $s\to \infty$.  In the case of ODE systems, the polynomial part is simply the constant feed-forward term, $\bD$.

 Let $\bH(s)$ be the transfer function of a DAE system.  We  decompose $\bH(s)$ as:
\begin{equation}
\bH(s) = \bC(s \bE - \bA)^{-1}\bB + \bD = \bG(s) + \bP(s),
\label{eq:GspP}
\end{equation}
where $\bG(s)$ is the strictly proper rational part, i.e., $\lim_{s\to \infty} \bG(s) = 0$ and $\bP(s)$ is the  polynomial part. Now, assume that 
the Petrov-Galerkin projection is applied to $\bH(s)$ as in (\ref{red_projection}).
Then, even though $\bE$ is singular, the reduced matrix $\bE_r = \bW ^T\bE \bV$ will generically be a~nonsingular matrix for $r \leq {\rm rank}(\bE)$. 
This means that unlike $\bH(s)$, which contains a polynomial part $\bP(s)$, the
reduced model will correspond to an ODE and the polynomial part of the reduced transfer function  $\bH_r(s)$ will be $\bD.$ Decompose $\bH_r(s)$ in a similar way 
$$
\bH_r(s) = \bG_r(s) + \bD,
$$
where $\bG_r = \bC_r(s \bE_r - \bA_r)^{-1}\bB_r$ is strictly proper. 
Then, the error transfer function 
$$
\bH(s) - \bH_r(s) = (\bG(s) - \bG_r(s)) + (\bP(s) - \bD)
$$
has a polynomial part $\bP(s) - \bD$.  Even when $\bP(s)$ is a polynomial of degree $1$, the difference $\bP(s) - \bD$ will grow without bound as $s \to \infty$, 
leading to unbounded $\hinf$ and $\htwo$ error norms. Even when $\bP(s)$ is a constant polynomial (i.e. degree $0$), unless $\bP(s) = \bD$,
this will still lead to unbounded $\htwo$ error. The only way to guarantee bounded error norms is to make sure that the reduced transfer function $\bH_r(s)$ has exactly the same polynomial part as $\bH(s)$, i.e. $\bH_r(s) = \bG_r(s) + \bP(s)$ so that the error function is simply $\bH(s) - \bH_r(s) = \bG(s) - \bG_r(s)$ having only a null polynomial component.  Based on these observations, \cite{GSW2013,wyatt2012issues} discusses how to modify the interpolatory projection bases $\bV$ and $\bW$
in order to achieve this goal.   As expected, the left and right deflating subspaces of the pencil $\lambda \bE - \bA$ corresponding to  finite and infinite eigenvalues play a crucial role:
\begin{thm} \label{interp_dae}
Suppose the transfer function $\bH(s) = \bC(s\bE-\bA)^{-1}\bB + \bD = \bG(s) + \bP(s)$ is associated with a DAE, where $\bG(s)$ and $\bP(s)$ are, respectively, the strictly proper and the polynomial parts of $\bH(s)$.
 Let  $\bbP_l$ 
and $\bbP_r$  be the spectral projectors onto the left and right deflating subspaces 
of the pencil $\lambda\bE-\bA$ corresponding to the finite eigenvalues. Also, let the columns 
of $\bW_\infty$ and $\bV_\infty$ span the left and right deflating subspaces of 
$\lambda\bE-\bA$ corresponding to the eigenvalue at infinity. 
For interpolation points $\sigma,\,\mu \in \IC$, suppose $\sigma\,\bE\, -\,\bA$ and $\mu\,\bE\, -\,\bA$ are invertible
and  $\rdir\in \IC^{m}$ and  $\ldir\in \IC^{\ell}$  are given (nontrivial)  tangent directions.
Suppose further that $\bH_r(s)$  is the reduced transfer function obtained by projection as in (\ref{red_projection}) using the model reduction bases $\bV$ and $\bW$. 
 Construct $\bV_{\!f}$ and 
$\bW_{\!f}$ so that 
\begin{align}
	\left(\left(\sigma\,\bE-\bA\right)^{-1}\bE \right)^{j-1}\left(\sigma\,\bE-\bA\right)^{-1} \bbP_l\bB\rdir\in  \mbox{\normalfont\textsf{Ran}}(\bV_{\!f}),
	 \label{eq:vf} 
	\end{align}
for $j=1,\ldots,N$,	and
	\begin{align}
\left(\left(\mu\,\bE-\bA\right)^{-T}\bE^T \right)^{j-1}\left(\mu\,\bE-\bA\right)^{-T} \bbP_r^T\bC^T\ldir\in  \mbox{\normalfont\textsf{Ran}}(\bW_{\!f})
 \label{eq:wf}
 \end{align}
 for $j=1,\ldots,M$.
  Define $\bW$ and $\bV$ using 
  $\bW =[\,\bW_{\!f}, \; \bW_\infty\,] \mbox{ and }  
\bV=[\,\bV_{\!f},\; \bV_\infty\,].$ Then, $\bH_r(s) = \bG_r(s) + \bP_r(s)$ satisfies $ {\bP_r}(s) = \bP(s)$ together with (\ref{eq:horti}) and (\ref{eq:holti}). If, in addition, $\sigma=\mu$, then 
(\ref{eq:hobti}) holds as well.
\end{thm}
Theorem \ref{interp_dae}, taken from \cite{GSW2013}, shows how to apply projection-based tangential interpolation to DAEs in the most general case where the index of the DAE and the interpolation points are arbitrary. By appropriately incorporating the deflating projectors $\bbP_r$ and $\bbP_l$ in the model reduction bases, the polynomial part of $\bH(s)$ is exactly matched as desired while simultaneously enforcing interpolation. For the  special case of DAEs  with proper transfer functions and  
interpolation around $s=\infty$, a solution has been given
in \cite{benner2006partial}. For descriptor systems of index~$1$, \cite{Ant2010imr}  offered a solution 
that uses an appropriately chosen $\bD_r$ term. 
\begin{rem}
A fundamental difficulty in the reduction of DAEs is the need to compute 
deflating projectors $\bbP_r$ and $\bbP_l$. 
For large-scale DAEs,  construction of 
$\bbP_r$ and $\bbP_l$ is at best very costly, if even feasible.  However, for the cases of semi-explicit descriptor systems of index-$1$ and Stokes-type descriptor systems of index~2,  it is possible
 to apply interpolatory projections without forming $\bbP_r$ and $\bbP_l$ explicitly  \cite{GSW2013,wyatt2012issues}.  Thus, no greater effort is required to produce reduced models than for the case of ODEs.  The Stokes-type descriptor systems of index~2 were also studied in \cite{heinkenschloss2008btm} in a balanced truncation setting.  Recently, \cite{ahmadbenner2013} extended the work of \cite{GSW2013} to address the reduction of index-$3$ 
DAEs without forming  projectors explicity.
For some structured problems arising in circuit simulation, multibody systems
or computational fluid dynamics, these projectors can be constructed without much computational effort \cite{Styk08}.
 We choose to omit these details from the present discussion but refer the interested reader to the original sources. 
\end{rem}
\section{Interpolatory Projections for $\htwo$ Optimal Approximation}  \label{sec:h2}
When interpolation points and tangent directions are specified, Section \ref{sec:intplt} presents an approach that one may follow in order to construct a reduced model  satisfying the desired  (tangential) conditions.   Notably, this development  does not suggest a strategy for choosing  interpolation points and tangent directions that lead to high-fidelity reduced models.  In this section, we approach this issue by developing interpolatory conditions that are necessary for  optimal approximation with respect to the $\htwo$ norm.
\subsection{Interpolatory $\htwo$-optimality Conditions} 
Consider the following optimization problem: Given a full-order system, $\bH(s)$, find a reduced model, $\bH_r(s)$ that minimizes the 
$\htwo$ error; i.e.,
\begin{equation} \label{h2opt} 
 \left\| \bH - \bH_r \right\|_{\htwo} =  \min_{ {\small
\dim(\widetilde{\bH}_{r})=r
}} \left\| \bH - \widetilde{\bH}_{r} \right\|_{\htwo}.
\end{equation}
As we pointed out in Section \ref{subsubsec:h2norm}, small  $\htwo$ error induces small time domain error $\left\| \by - \by_r \right\|_{L_\infty}$,
so attempting to \emph{minimize} the $\htwo$ error is a worthy goal. 

The $\htwo$ optimization problem (\ref{h2opt}) is nonconvex ; finding a global minimizer
is typically infeasible. A common approach used instead involves finding {\em locally} optimal reduced models that satisfy first-order necessary conditions for optimality.  The (simpler) problem of finding locally optimal $\htwo$ reduced models has been studied extensively.   Optimality conditions have been formulated either in terms of Lyapunov and Sylvester equations \cite{wilson1970optimum,hyland1985theoptimal,halevi1992fwm,zigic1993contragredient,yan1999anapproximate,spanos1992anewalgorithm,benner_sylvester} or in terms of rational (tangential) interpolation conditions \cite{meieriii1967approximation,gugercin2005irk,gugercin2006rki, gugercin2008hmr, vandooren2008hom,bunse-gerstner2009hom, kubalinska2007h0i, beattie2007kbm,beattie2009trm,beattie2012realization,krajewski1995program,panzer2013}.   \cite{gugercin2008hmr} showed the equivalence between the Lyapunov/Sylvester equation conditions and the interpolation framework that we describe here.

We will first assume that $\bE$ is nonsingular, so that $\bH(s) = \bC(s \bE -\bA)^{-1}\bB+\bD $ will correspond to a system of ODEs and
 $\lim_{s\to \infty} \bH(s) = \bD$. In order to have a bounded $\htwo$ error norm, $\| \bH - \bH_r\|_\htwo$,  it is necessary that $\bD_r = \bD$. Therefore, without loss of generality, one may take $\bD_r = \bD = \mathbf{0}$.

For MIMO systems, interpolatory first-order conditions for $\htwo$ optimality are best understood 
from the pole-residue expansion for $\bH_r(s)$. We  write $\bH_r(s) = \bC_r(s \bE_r -\bA_r)^{-1}\bB_r$ in the following way:
\begin{equation} \label{eq:pr}
\bH_r(s)=\sum_{i=1}^r\frac{\ldir_i\rdir_i^T}{s-{\lambda}_i}
\end{equation}
where we have assumed that the $\lambda_i$s are distinct. We refer to $\ldir_i \in \IC^{p}$ and $\rdir_i \in \IC^{m}$ in (\ref{eq:pr}) respectively as left/right residue directions associated with the pole ${\lambda}_i$ of $\bH_r(s)$; $\ldir_i\rdir_i^T$ is the (matrix) residue of $\bH_r(s)$ at $s=\lambda_i$. The pole-residue expansion in (\ref{eq:pr}) can be computed effectively by computing a generalized eigenvalue decomposition for the matrix pencil $\lambda \bE_r - \bA_r$, which is a trivial computation for the small to modest orders of $r$ typically encountered. Note that finding such a representation for the full model $\bH(s)$ will generally be infeasible.
\begin{thm}  \label{h2cond}
Let $\bH_r(s)$ in (\ref{eq:pr}) be the best $r{\rm th}$ order rational
 approximation of $\bH(s)$ with respect to the $\htwo$ norm.
Then, 
\begin{subequations}\label{H2optcond}
\begin{align}
 \bH(-{\lambda}_k) \rdir_k &=\bH_r(-{\lambda}_k) \rdir_k,\label{H2optcondright} \\
 \ldir_k^T \bH(-{\lambda}_k) &=  \ldir_k^T\bH_r(-{\lambda}_k),  \mbox{ and }
 \label{H2optcondleft} \\
 \ldir_k^T \bH'(-{\lambda}_k) \rdir_k 
   &=  \ldir_k^T\bH_r'(-{\lambda}_k) \rdir_k \label{H2optcondhermite} 
\end{align}
\end{subequations}
for $k=1,\,2,\,...,\,r$.
\end{thm}
In particular, any optimal $\htwo$ approximation $\bH_r(s)$ must be a bitangential Hermite interpolant to $\bH(s)$, and this theorem directly connects optimal model reduction to interpolation. The optimal interpolation points and tangent directions are derived from the pole-residue representation of $\bH_r(s)$: The optimal interpolation points are the mirror images of the poles of $\bH_r(s)$ reflected across the imaginary axis, and the optimal tangent directions are the residue directions associated with that pole. 

Interpolatory conditions for SISO systems were initially introduced by Meier and Luenberger \cite{meieriii1967approximation}.
However, until recently effective numerical algorithms to find reduced systems that satisfy these conditions were lacking, especially for large-scale settings.  Gugercin {\it et al.} in \cite{gugercin2005irk, gugercin2006rki} 
introduced such an algorithm, called the \textit{Iterative Rational Krylov Algorithm} (\IRKA). In practice, \IRKA~has significantly expanded the utility of  optimal $\htwo$ model reduction.  The optimality conditions for MIMO systems as presented in Theorem \ref{h2cond} were developed in \cite{gugercin2008hmr,bunse-gerstner2009hom,vandooren2008hom}, and led to an analogous algorithm for \IRKA~in the MIMO case.  This is the main focus of Section \ref{sec:irka}.  Recall that we have assumed that $\bH_r(s)$ has distinct (reduced) poles: $\lambda_1,\ldots,\lambda_r$.   Optimality conditions for cases when $\bH_r(s)$ has repeated poles are derived in 
\cite{van2010h_2}.

\subsection{IRKA for optimal $\htwo$ approximation} \label{sec:irka}
Theorem \ref{h2cond} gives optimality conditions that depend on the poles and residues of a reduced order system,
a locally $\htwo$-optimal reduced system, that is unknown {\em a priori } and is yet to be determined.    \IRKA\ utilizes the construction of Theorem \ref{thm:interpolation} to force interpolation at the mirror images of successive sets of reduced poles, iteratively correcting the reduced model until the optimality conditions of Theorem \ref{h2cond} hold. 
The method proceeds as follows: Given some initial interpolation points $\{\sigma_i\}_{i=1}^r$ and directions
$\{\rdir_i\}_{i=1}^r$ and $\{\ldir_i\}_{i=1}^r$, construct $\bV$ and $\bW$ as in (\ref{eqn:Vr}) and (\ref{eqn:Wr}), respectively, and construct an intermediate reduced model $\bH_r(s)$ using (\ref{red_projection}). Then, compute the pole-residue decomposition of  $\bH_r(s)$, 
 $$\bH_r(s) = \sum_{i=1}^r\frac{\widehat{\ldir}_i\widehat{\rdir}_i^T}{s-{\lambda}_i}$$
 (by solving a small $r \times r$ generalized eigenvalue problem). For $\bH_r(s)$ to satisfy the first-order necessary conditions, we need $\sigma_i = -\lambda_i$, $\rdir_i = \widehat{\rdir}_i$, and 
 $\ldir_i = \widehat{\ldir}_i$, for
$i=1,\ldots,r$.
 Therefore, set 
$$
\sigma_i \longleftarrow -\lambda_i,\quad\rdir_i \longleftarrow \widehat{\rdir}_i,\quad\mbox{ and } 
\quad\ldir_i \longleftarrow \widehat{\ldir}_i,\quad \mbox{for}\quad i=1,\ldots,r
$$
as the next interpolation data and iterate until convergence is reached. A brief sketch of \IRKA~is given below:
\begin{center}
    \framebox[4.5in][t]{
\vspace*{-0.5cm}
    \begin{minipage}[c]{4.3in}
    {\small
    \begin{alg} \label{sucratkry} 
    {\bf MIMO $\htwo$-Optimal Tangential Interpolation (``\IRKA")}        
    \begin{enumerate}
       \item Make an initial $r$-fold shift selection: $\{\sigma_1,\ldots,\sigma_r\}$ that is closed under conjugation      
       (i.e., $\{\sigma_1,\ldots,\sigma_r\}\equiv\{\overline{\sigma_1},\ldots, \overline{\sigma_r}\}$ viewed as sets) 
        and initial \\ tangent directions 
       ${\rdir}_1,\ldots,{\rdir}_r$ and 
       ${\ldir}_1,\ldots,{\ldir}_r$, also closed under conjugation.
       \item  $\bV_r =  \left[(\sigma_1 \bE- \bA)^{-1}\bB{\rdir}_1~\cdots~(\sigma_r \bE - \bA)^{-1}\bB{\rdir}_r~\right]$
      \item  $\bW_r =  \left[({\sigma_1}\, \bE - \bA^T)^{-1}\bC^T{\ldir}_1 \,\cdots \,({\sigma_r}\, \bE - \bA^T)^{-1}\bC^T{\ldir}_1~\right]$
       \item while (not converged)
       \begin{enumerate}
       \item $\bA_r = \bW_r^T \bA \bV_r$, $\bE_r = \bW_r^T \bE \bV_r$, $\bB_r = \bW_r^T\bB$, and $\bC_r = \bC \bV_r$
       \item Compute a pole-residue expansion of $\bH_r(s)$:
         $$
         \bH_r(s) = \bC_r(s \bE_r -\bA_r)^{-1}\bB_r = \sum_{i=1}^r
         \frac{\widehat{\ldir}_i \widehat{\rdir}_i^T}{s-{\lambda}_i}
         $$
              \item $\sigma_i \longleftarrow -\lambda_i,\quad\rdir_i \longleftarrow \widehat{\rdir}_i,\quad\mbox{ and } 
\quad\ldir_i \longleftarrow \widehat{\ldir}_i,\quad \mbox{for}\quad i=1,\ldots,r$
     \item  $\bV_r =  \left[(\sigma_1 \bE- \bA)^{-1}\bB{\rdir}_1~\cdots~(\sigma_r \bE - \bA)^{-1}\bB{\rdir}_r~\right]$
      \item  $\bW_r =  \left[({\sigma_1}\, \bE - \bA^T)^{-1}\bC^T{\ldir}_1 \,\cdots \,({\sigma_r}\, \bE - \bA^T)^{-1}\bC^T{\ldir}_1~\right]$
       \end{enumerate}
       \item $\bA_r = \bW_r^T \bA \bV_r$, $\bE_r = \bW_r^T \bE \bV_r$, $\bB_r = \bW_r^T \bB$, $\bC_r = \bC\bV_r$
       \end{enumerate}
       \end{alg}
       }
    \end{minipage}
    }
    \end{center}
  Upon convergence, the reduced model, $\bH_r(s)$, satisfies the interpolatory first-order necessary conditions (\ref{H2optcond}) for 
  $\htwo$ optimality by construction.  Convergence is generally observed to be rapid; though it slows as input/output orders grow. Convergence may be guaranteed {\em a priori} in some circumstances \cite{flagg2012convergence}; yet there are known cases where convergence may fail \cite{gugercin2008hmr,flagg2012convergence}.  When convergence occurs, the resulting reduced model is guaranteed to be a local $\htwo$-minimizer since the local maxima of the $\htwo$ minimization problem are known to be repellent \cite{krajewski1995program}.  Overall in practice,
 \IRKA\ has seen significant success in computing high fidelity (locally) optimal reduced models and 
 has been successfully applied in large-scale settings to find $\htwo$-optimal reduced models for systems with hundreds of thousands of state variables; for example, see \cite{KRXC08} for application in  cellular neurophysiology, \cite{borggaard2012model} for energy efficient building design in order to produce accurate compact models for the indoor-air environment, \cite {gugercin2008hmr} for optimal cooling for steel profiles. Moreover, \cite{Breiten_H2} has extended \IRKA~to the reduction of \emph{bilinear} dynamical systems,  a special class of weakly nonlinear dynamical systems.

Our analysis so far has assumed that $\bE$ is a nonsingular matrix. Interpolatory optimal $\htwo$ model reduction for the case of singular $\bE$, i.e., for systems of DAEs, has been developed in \cite{GSW2013} and \IRKA~has been extended to DAEs. Similar to the ODE case where we require $\bD = \bD_r$, the DAE case requires that the polynomial part of $\bH_r(s)$ match that of $\bH(s)$ exactly and the strictly proper part of $\bH_r(s)$ be an optimal tangential interpolant to the strictly proper part of $\bH(s)$. For details, we refer the reader to \cite{GSW2013}.

\subsection{Interpolatory Weighted-$\htwo$ Model Reduction}
The error measures we have considered thus far give the same weight to all frequencies equally and they are global in nature to the extent that
  degradation in fidelity is penalized in the same way throughout the full frequency spectrum.  However, some applications require that  certain frequencies be weighted more than others. For example, certain dynamical systems, such as mechanical systems or electrical circuits,  might operate only in certain frequency bands and retaining fidelity outside this frequency band carries no value.  
  This leads to the problem of weighted model reduction. We formulate it here in terms of a weighted-$\htwo$ norm. 

Let $\bW(s)$ be an input weighting function, a ``shaping filter.''  
We will assume that $\bW(s)$ is a rational function itself in the form 
 \begin{align}
\bW(s)=\mathbf{C}_w\left(s\mathbf{I}-
\mathbf{A}_w\right)^{-1}\mathbf{B}_w+\bD_w =\sum_{k=1}^{n_w} \frac{\be_k\ \boldf_k^T}{s-\gamma_k} +
\bD_w, \label{WPoleRes}
 \end{align}
where $n_w$ denotes the dimension of $\bW(s)$.  Then, given the full-model $\bH(s)$ and the weight $\bW(s)$, define the weighted $\htwo$ error norm as:
\begin{equation}\label{wtdError}
\| \bH- \bH_r \|_{{\htwo}(W)} \stackrel{\tiny{\mbox{def}}}{=} \| \left( \bH - \bH_r \right)\cdot
\bW\|_{\htwo} .
\end{equation}
In addition to the input weighting $\bW(s)$, one may also define a filter for output weighting.
 For simplicity of presentation, we focus here on one-sided weighting only.  The goal is to find a 
reduced model $\bH_r(s)$ that minimizes the weighted error (\ref{wtdError}):
\begin{equation} \label{eq:H2optProb}
\| \bH - \bH_r \|_{{\htwo}(W)}  =
\min_{ {\small\dim(\widetilde{\bH}_{r})=r}} \|\bH-\widetilde{\bH}_r\|_{{\htwo}(W)}. 
\end{equation}
Weighted-$\htwo$ model reduction has been considered in \cite{halevi1992fwm} and  \cite{spanos1992anewalgorithm} 
using a framework that uses Riccati and/or Lyapunov equations.  A numerically more efficient, interpolation-based approach was  introduced in \cite{Ani13} for the SISO case.  This initial interpolatory framework was significantly extended and placed on a more rigorous theoretical footing (which allowed straightforward extension to MIMO systems) in \cite{breiten2013near} where the equivalence of the Riccati and interpolation-based frameworks were also proved. Our presentation below follows \cite{breiten2013near}. We use the notation $\bH \in \htwo$ to indicate that $\bH(s)$ is a stable dynamical system with $\| \bH \|_\htwo < \infty$. We use the notation $\bH \in \htwo(W)$  analogously.

The interpolatory framework for the weighted-$\htwo$ problem is best understood by defining a new linear transformation (\cite{breiten2013near,Ani13})
 { \begin{align} 
\mathfrak{F}[\bH](s)= & \bH(s)\bW(s)\bW(-s)^T \label{Fmap}  
+ \sum^{n_w}_{k=1}
\bH(-\gamma_{k})\bW(-\gamma_{k})\frac{\boldf_k\be_k^T}{s+\gamma_{k}} ,
\end{align}
}
where $\bH\in \htwo(W)$, and $\boldf_k$ and $\be_k$ are as defined in (\ref{WPoleRes}).  
$\mathfrak{F}[\bH](s)$ is a bounded linear transformation from $\htwo(W)$ to $\htwo$
\cite{breiten2013near}.  A state-space representation for $\mathfrak{F}[\bH](s)$ is given by

{ 
\begin{align}\label{eq:trans_func_repr}
  \mathfrak{F}[\bH](s) &= \cbC_\mathfrak{F} (s \bI - \cbA_\mathfrak{F})^{-1} \cbB_\mathfrak{F}\\
 &=\underbrace{\begin{bmatrix} \bC & \bD\bC_w \end{bmatrix}}_{\cbC_{\mathfrak{F}}}
 \Bigg( s\bI - 
 \underbrace{\begin{bmatrix} \bA & \bB\bC_w \\ \mathbf{0} & \bA_w
\end{bmatrix}}_{{\cbA_\mathfrak{F}}}\Bigg)^{-1}
\underbrace{ \begin{bmatrix}\bZ\bC_w^T + \bB\bD_w\bD_w^T \\ \bP_w \bC_w^T +
\bB_w\bD_w^T \end{bmatrix}}_{\cbB_{\mathfrak{F}}}, \nonumber
\end{align}}
where $\bP_w$ and $\bZ$  solve, respectively, 
\begin{align}
  \bA_w & \bP_w + \bP_w \bA_w^T + \bB_w\bB_w^T = \mathbf{0} \quad \mbox{and} \label{eq:weight_lyap} \\
  \bA   \bZ& + \bZ \bA_w^T + \bB (\bC_w\bP_w+\bD_w\bB_w^T) =\mathbf{0}. \label{eq:weight_sylv}
\end{align}
For $\bH,\ \bH_r\in \htwo(W)$, 
 denote the impulse responses corresponding to    
$\mathfrak{F}[\bH](s)$ and $\mathfrak{F}[\bH_r](s)$, respectively by $\bF(t)$ and $\bF_r(t)$, so that 
$\mathfrak{F}[\bH]=\mathcal{L}\left\{\bF\right\}$ and 
$\mathfrak{F}[\bH_r]=\mathcal{L}\left\{\bF_r\right\}$, where 
$\mathcal{L}\left\{\cdot\right\}$ denotes the Laplace transform. 

\begin{thm} \label{thm:weight_intp_cond}
Given the input weighting $\bW(s)$, let $\bH_r\in \htwo(W)$ be the best order-$r$ rational approximation to $\bH$ in the weighted-$\htwo$ norm. 
Suppose that $\bH_r$ has the form 
{
\begin{equation}
\bH_r(s)=\mathbf{C}_r\left(s\mathbf{I}-
\mathbf{A}_r\right)^{-1}\mathbf{B}_r +\bD_r
 =\sum_{k=1}^{n_r} \frac{\ldir_k\ \rdir_k^T}{s-\lambda_k} +\bD_r. \label{GrRepr}
 \end{equation}}
 where $\lambda_1,\ldots,\lambda_r$ are assumed to be distinct.
Then 
\begin{subequations}\label{eq:weight_intp_cond}
\begin{align}
 \mathfrak{F}[\bH](-\lambda_k)\rdir_k &= \mathfrak{F}[\bH_r](-\lambda_k) \rdir_k
\label{rightInterpCond}\\[.1in]
\ldir_k^T\mathfrak{F}[\bH](-\lambda_k) &=\ldir_k^T \mathfrak{F}[\bH_r](-\lambda_k), \mbox{ and }
\label{leftInterpCond} \\[.1in]
\ldir_k^T\mathfrak{F}^{\,\prime}[\bH](-\lambda_k)\rdir_k &=
\ldir_k^T\mathfrak{F}^{\,\prime}[\bH_r](-\lambda_k) \rdir_k.    \label{bitanInterpCond}  \\[.1in]
 \label{IRFInterpCond}
\bF(0)\bn&=\bF_r(0)\bn.
 \end{align}
  \end{subequations}
  for $k=1,\,2,\,...,\,r$  and for all $\bn\in \mathsf{Ker}(\bD_w^T)$
where 
$\mathfrak{F}^{\,\prime}[\,\cdot\,](s) = \frac{d\
}{ds}\mathfrak{F}[\,\cdot\,](s)$.
\end{thm}
Bitangential Hermite interpolation once again appears as a necessary condition for optimality. However, unlike the unweighted $\htwo$ case, the interpolation conditions need to be satisfied by the maps 
$\mathfrak{F}[\bH](s)$ and  $\mathfrak{F}[\bH_r](s)$ as opposed to $\bH(s)$ and $\bH_r(s)$. For 
$\bW(s) = \bI$,  (\ref{rightInterpCond})-(\ref{bitanInterpCond}) simplify to (\ref{H2optcondright})-(\ref{H2optcondhermite}) and  (\ref{IRFInterpCond}) is automatically satisfied since $\mathsf{Ker}(\bD_w^T) = \{\mathbf{0}\}$ when $\bW(s) = \bI$.

Guided by how \IRKA~is employed to satisfy the $\htwo$ optimality conditions (\ref{H2optcondright})-(\ref{H2optcondhermite}), one might consider simply applying \IRKA~to the state-space representation of $\mathfrak{F}[\bH](s)$ given in (\ref{eq:trans_func_repr}) to satisfy the weighted interpolation conditons (\ref{rightInterpCond})-(\ref{bitanInterpCond}). Unfortunately, this solves a different problem: If \IRKA~is applied to $\mathfrak{F}[\bH](s)$ directly, then one obtains a reduced model $\bH_r(s)$ that interpolates $\mathfrak{F}[\bH](s)$. That is, instead of 
(\ref{rightInterpCond}), one obtains instead  $\mathfrak{F}[\bH](-\lambda_k)\rdir_k = \bH_r(-\lambda_k) \rdir_k$ which is clearly not appropriate. The need to preserve the structure in the maps $\mathfrak{F}[\bH](s)$ and  $\mathfrak{F}[\bH_r](s)$ while satisfying  interpolatory conditions makes
the development of an \IRKA-like algorithm for the weighted-$\htwo$ problem quite nontrivial. 
Breiten {\it et al} in \cite{breiten2013near} proposed an algorithm, called  ``Nearly Optimal Weighted Interpolation" (\NOWI), that \emph{nearly satisfies} the interpolatory optimality conditions 
 (\ref{rightInterpCond})-(\ref{bitanInterpCond}) while preserving the structure in $\mathfrak{F}[\bH](s)$ and  $\mathfrak{F}[\bH_r](s)$.  The deviation from exact interpolation is quantified explicitly. 
Even though \NOWI~proves itself to be a very effective numerical algorithm in many circumstances, the development of an algorithm that satisfies (\ref{rightInterpCond})-(\ref{bitanInterpCond}) exactly remains an important future goal.

\subsection{Descent Algorithms for $\htwo$ Model Reduction}
At its core, \IRKA~is a fixed point iteration.  Excepting the  special case of symmetric state-space systems,\footnote{$\bE = \bE^T$ is positive definite,  $\bA = \bA^T$
is negative definite, and $\bB = \bC^T$} where convergence is guaranteed,   
convergence of \IRKA~for general systems is not guaranteed, 
see \cite{gugercin2008hmr,flagg2012convergence}, although superior performance is commonly observed.
More significantly, \IRKA~is not a descent algorithm; that is, the $\htwo$ error 
might fluctuate during intermediate steps and premature termination of the algorithm could result (at least in principle) in a worse approximation than what was provided for initialization. 
 To address these issues, Beattie and Gugercin \cite{beattie2009trm} developed an $\htwo$ descent algorithm that reduces the $\htwo$ error  at each step of the iteration and assures global convergence to a local minimum. 

The key to their approach is the following representation of the $\htwo$ error norm for MIMO systems  \cite{beattie2009trm}:
\begin{thm} \label{h2_exact}
Given a full-order model $\bH(s)$, let $\bH_r(s)$  have the form in (\ref{eq:pr}), i.e.,
$$\bH_r(s) = \sum_{i=1}^r \frac{{\ldir}_i{\rdir_i}^T}{s-{\lambda}_i}.$$
Then, the $\htwo$
norm of the error system  is given by
\begin{align} 
\left\| \bH-\bH_r \right\|_{\htwo}^2 = \left\| \bH \right\|_{\htwo}^2  - 2 
\sum_{k=1}^r \,{\ldir}_k^T\bH(-{\lambda}_k){\rdir}_k  \label{h2normexp}
 + \sum_{k,\ell=1}^r
\frac{{\ldir}_k^T{\ldir}_\ell\,{\rdir}_\ell^T{\rdir}_k}
{-{\lambda}_k-{\lambda}_\ell}.
\end{align}
\end{thm}
For SISO systems, Krajewski {\it et al.} \cite{krajewski1995program} developed and proved
a similar expression, which was later rediscovered in \cite{gugercin2003anmathcal,gugercinprojection,antoulas2005approximation}.

If one considers  $\bH_r(s) = \sum_{i=1}^r \frac{{\ldir}_i{\rdir_i}^T}{s-{\lambda}_i}$ in the pole-residue form, the variables defining the reduced model are the residue directions ${\ldir}_i$, ${\rdir_i}$ and the poles ${\lambda}_i$, for $i=1,\ldots,r$.  The formula (\ref{h2normexp})  expresses the error in terms of these variables. Thus, one can compute the gradient and Hessian of the error with respect to unknowns and construct globally convergent descent (optimization) algorithms. Gradient and Hessian expressions were derived in \cite{beattie2009trm}. For brevity of presentation, we  include only the gradient expressions here.
\begin{thm} \label{gradJ}
Given the full-order model $\bH(s)$ and the reduced model $\bH_r(s)$ as in (\ref{eq:pr}), define 
$$ \cJ \stackrel{\tiny{\mbox{\rm def}}}{=}  \left\| \bH-\bH_r \right\|_{\htwo}^2. $$
Then, for $i=1,\ldots,r$,
  \begin{align}
\frac{\partial \cJ}{\partial \lambda_i}&\,=\, - 2\,\ldir_i^T
\left(\bH_r'(-\lambda_i) 
-\bH'(-\lambda_i) \right)\rdir_i  . \label{djdlambda}
\end{align}
Moreover, the gradient of $\cJ$ with respect to residue directions listed as
$$
\{\rdir,\ldir\}
=[\rdir_1^T, \ldir_1^T, \rdir_2^T, \ldir_2^T,\ldots,\rdir_r^T, \ldir_r^T]^T,
$$
 is given by $\nabla_ {\!\!\{\rdir,\ldir\} }\cJ$, a vector of length $r(m+p)$, partitioned into $r$ vectors of length $m+p$ as
 \begin{align}
\left(\nabla_{\!\!\{\rdir,\ldir\} }\cJ\right)_\ell=\left(
\begin{array}{c}
2\,  \left(\ldir_\ell^T\bH_r(-\lambda_\ell)-\ldir_\ell^T\bH(-\lambda_\ell)\right)^T\\
2\,  \left(\bH_r(-\lambda_\ell) \rdir_\ell-\bH(-\lambda_\ell) \rdir_\ell \right)
\end{array}\right)  \label{dJdbc} 
\end{align}
for $\ell=1,\,2,\,\ldots,\,r.$  
\end{thm}
One may observe that setting the gradient expression in  (\ref{djdlambda}) and (\ref{dJdbc}) to zero leads immediately to the interpolatory optimality conditions (\ref{H2optcond}). Having gradient and Hessian expressions at hand, one may then develop a globally convergent descent algorithm for $\htwo$ optimal reduction as done in \cite{beattie2009trm} where the optimization algorithm was put in a trust-region framework. Unlike in \IRKA, the intermediate reduced models are not interpolatory. However, upon convergence, they satisfy the interpolatory optimality conditions. 

In a recent paper, for SISO systems $H(s) = \bc^T (s \bE - \bA)^{-1}\bb$ where 
$\bb,\bc \in \IR^n$ are length-$n$ vectors, Panzer {\it et al} \cite{panzer2013} applied a descent-type algorithm  successively. Instead of designing a degree-$r$ rational function directly,  \cite{panzer2013}  first constructs a SISO degree-$2$ rational function $H_r(s) = \bc_r^T (s \bE_r - \bA_r)^{-1}\bb_r$  where $\bA_r, \bE_r \in \IR^{2 \times 2}$ and $\bb_r,\bc_r^T \in \IR^2$ by a descent method where only  Lagrange optimality conditions are enforced (without Hermite conditions). 
Then, the error transfer function is decomposed in a multiplicative form
$$
H - H_r(s) = \underbrace{\left( \bc^T (s \bE - \bA)^{-1}\bb_\perp \right)}_{H_\perp(s)} \left({\bc}_r^T (s \bE_r - \bA_r)^{-1}\bb_r \right)
$$
where $\bb_\perp = \bb - \bE \bV (\bW^T \bE \bV)^{-1}\bb_r$ and the method proceeds by constructing another degree-$2$ approximation to $H_\perp(s)$ in a descent framework, once more only enforcing the Lagrange optimality conditions. At the end, all the intermediate degree-$2$ approximants are put together in a special way to form the final reduced model of degree-$r$. For details, we refer the reader to  \cite{panzer2013}.  The final reduced model will not generally satisfy the full set of interpolatory $\htwo$ optimality  conditions - only the Lagrange conditions are satisfied.  Moreover, the incremental approach means optimization over a smaller set; 
thus for a given $r$, optimization directly over a degree-$r$ rational function using 
the gradient and Hessian expressions in \cite{beattie2009trm}  will lead to a smaller model reduction error than an incremental search. However since \cite{panzer2013} performs the optimization over a very small number of variables in each step, this approach can provide some numerical advantages. 

Druskin {\it et al} in \cite{druskin2011adaptive} and \cite{druskin2012adaptive} suggest alternative greedy-type algorithms for interpolatory model reduction.  Instead of constructing $r$ interpolation points (and directions at every step) as done in \IRKA~or in the descent framework of \cite{beattie2009trm},
 \cite{druskin2011adaptive} and \cite{druskin2012adaptive} start instead with an interpolation point and corresponding tangent directions. Then, a greedy search on the residual determines the next set of interpolation data. Since the greedy search is not done in a descent setting, this is not a descent method and at the end optimality conditions will not be satisfied typically. Nonetheless, the final reduced model is still an interpolatory method. Even though the resulting reduced models will not generally be as accurate as those obtained by \IRKA,~the methods of  \cite{druskin2011adaptive} and \cite{druskin2012adaptive}  provide satisfactory approximants at relatively low cost. 
 
 Descent-type algorithms have been extended to the weighted-$\htwo$ norm minimization as well; for details see \cite{petersson2013nonlinear,vuillemin2012spectral,Bre13}.

\section{Interpolatory Model Reduction for Generalized Coprime Framework}  \label{sec:coprime}
So far, we have assumed that the original transfer function has a generic first-order state-space representation: $\bH(s) = \bC(s\bE - \bA)^{-1}\bB$. This representation is quite general and a wide range of linear dynamical systems can be converted to this form, at least in principle.
However, problem formulations often lead to somewhat different structures that reflect the underlying physics or other important system features.  
One may wish to retain such structural features and conversion to standard first-order representations often obfuscates these features and
may lead then to ``unphysical" reduced models.   Neutral and delay differential equations present another interesting class of dynamical systems that are generically of infinite order, so they do not accept a standard first-order representation using a finite dimensional state space.  In this section, we 
follow the discussion of Beattie and Gugercin \cite{beattie2008ipm}, and show how interpolatory methods can be used to preserve relevant system structure in reduced models, often avoiding entirely the need to convert the system to an equivalent first-order state-space representation.

\subsection{Polynomial and Delay Systems }
A common example of a situation where conversion to the standard state-space form is 
possible but may not be prudent is the case of constant coefficient ordinary differential equations of order two or more, with dynamics given by 
\begin{align} \label{fompoly}
\bA_{0}\frac{d^{\ell}\bx}{dt^{\ell}}  &+ \bA_{1}\frac{d^{\ell-1}\bx}{dt^{\ell-1}} +\cdots
 + \bA_{\ell}\bx(t) = \bB \bu(t), \\ 
 \by(t)&= \bC_{1}\frac{d^{\ell-1}\bx}{dt^{\ell-1}}  + \bC_{2}\frac{d^{\ell-2}\bx}{dt^{\ell-2}} +\cdots
 + \bC_{\ell}\bx(t),  \nonumber
\end{align}
where $\bA_i \in \IR^{n \times n}$, for $i=0,\ldots,\ell$, $\bB \in \IR^{n\times m}$ and $\bC_i \in \IR^{p \times n}$  for $i=1,\ldots,\ell$. By defining a state vector $\bq = [\bx^T,~\dot{\bx}^T,\ddot{\bx}^T,\ldots,(\bx^{(\ell-1)})^T]^T$, one may easily convert (\ref{fompoly}) into an equivalent first-order form (\ref{ltisystemintro}). However, this has two major disadvantages: 
\begin{enumerate}
\item By forming the vector $\bq(t)$, the physical meaning of the state vector $\bx(t)$ is lost in the model reduction state since the model reduction process will mix the physical quantities such as displacement and velocity; the reduced state loses its physical significance. 
\item  The dimension of the aggregate state $\bq(t)$ is  $\ell \times n$. As a consequence, conversion to first-order form has made the model reduction problem numerically much harder. For example, if reduction is approached via interpolation, we now need to solve linear systems of size 
$(\ell\,n) \times (\ell\,n)$
\end{enumerate}
Therefore, it is desirable to perform model reduction in the original state-space associated with the original  representation (\ref{fompoly});
we wish to preserve the structure of (\ref{fompoly}) in the reduced model and produce a  reduced model of the form 
\begin{align} \label{fompolyrom}
\bA_{0,r}\frac{d^{\ell}\bx_r}{dt^{\ell}}  &+ \bA_{1,r}\frac{d^{\ell-1}\bx_r}{dt^{\ell-1}} +\cdots
 + \bA_{\ell,r}\bx_r(t) = \bB_r \bu(t)  \\
  \by(t)&= \bC_{1,r}\frac{d^{\ell-1}\bx_r}{dt^{\ell-1}}  + \bC_{2,r}\frac{d^{\ell-2}\bx_r}{dt^{\ell-2}} +\cdots
 + \bC_{\ell,r}\bx_r(t)  \nonumber
 \end{align}
where $\bA_{i,r} \in \IR^{r \times r}$, for $i=0,\ldots,\ell$, $\bB_r \in \IR^{r\times m}$ and $\bC_{i,r} \in \IR^{p \times r}$  for $i=1,\ldots,\ell$. 

Another example where the structure of a dynamical system presents an obstacle to reduction using methods that depend on availability of a standard first-order form is the class of delay differential equations. Consider a linear dynamical system with an internal delay, given in state space
form as:
\begin{equation} \label{ExDelaySS}
\bE \dot{\bx}(t)  = \bA_{0}\, \bx(t) + \bA_{1}\, \bx(t-\tau) 
+ \bB\, \bu(t),  
 \quad \by(t)  =  \bC\, \bx(t) 
\end{equation}
with $\tau>0$, and $\bE,\,\bA_{0},\, \bA_{1}  \in \IR^{n \times n}$, $\bB \in \IR^{n\times m}$ and 
$\bC \in \IR^{p \times n}$. The system in (\ref{ExDelaySS}) is not associated with a rational transfer function due to the delay term; it is intrinsically of infinite order. Preserving the delay structure in the reduced model is crucial for accurate representation and so, we seek a reduced model of the form 
\begin{equation} \label{ExDelaySSrom}
\bE_r \dot{\bx}_r(t)  = \bA_{0_r}\, \bx_r(t) + \bA_{1,r}\, \bx_r(t-\tau) 
+ \bB_r\, \bu(t),  
 \quad \by_r(t)  =  \bC_r\, \bx_r(t) 
\end{equation}
with $\tau>0$, and $\bE_r,\,\bA_{0,r},\, \bA_{1,r}  \in \IR^{r \times r}$, $\bB_r \in \IR^{r\times m}$ and 
$\bC \in \IR^{p \times r}$. We want to perform this reduction step without the need for approximating the delay term with an additional rational approximation.

\subsection{Generalized Coprime Representation}
Examples such as these lead us to consider transfer functions having the following \emph{Generalized Coprime Representation}: 
\begin{equation}    \label{HdecompD}
     \cbH(s)=\cbC(s)\cbK(s)^{-1}\cbB(s) 
     + \cbD
\end{equation}
where $\cbD$ 
is a constant $p \times m$ matrix, both $\cbC(s)\in \IC^{p \times n}$ 
and $\cbB(s)\in \IC^{n\times m}$ 
are analytic in the right half plane, and  
$\cbK(s)\in \IC^{n \times n}$ is analytic and full rank 
throughout the right halfplane.  Note that both (\ref{fompoly}) and (\ref{ExDelaySS}) fit this framework:
for the polynomial system (\ref{fompolyrom}), we obtain
%\begin{equation} \label{eq:polyH}
%\cbH(s)= \underbrace{s^{\ell-1}\bC_1+s^{\ell-2}\bC_2 + \cdots + s \bC_{\ell-1} + \bC_\ell)}_{\cbC(s)} \underbrace{(s^\ell \bA_0 + s^{\ell -1} \bA_1 + \cdots + s \bA_{\ell-1} + \bA_\ell)}_{\cbK(s)}\,^{-1} \underbrace{\bB}_{\cbB(s)},
%\end{equation}
\begin{equation} \label{eq:polyH}
\cbH(s)= \underbrace{\left(\sum_{i=1}^\ell s^{\ell-i}\bC_i\right)}_{\cbC(s)} \underbrace{
\left(\sum_{i=0}^\ell s^{\ell-i}\bA_i\right)^{-1}}_{{\displaystyle \cbK(s)^{-1}}} \,\,\underbrace{\bB}_{\cbB(s)},
\end{equation}
and for the delay system  (\ref{ExDelaySS}), we have
\begin{equation} \label{eq:delayH}
 \cbH(s) = \underbrace{\bC}_{\cbC(s)} \underbrace{\left(s\,\bE -\bA_{0} -e^{-\tau s}\, \bA_{1}\right)^{-1}}_{{\displaystyle \cbK(s)^{-1}}}\, \underbrace{\bB}_{\cbB(s)}.
\end{equation}
Our model reduction goals are the same: Construct a reduced (generalized coprime) transfer function that tangentially interpolates the original one. To this end, we choose two model reduction bases  $\cbV \in {\mathbb R}^{n \times r}$ and $\cbW\in {\mathbb R}^{n \times r}$ as before. This  leads to a reduced transfer function 
\begin{equation}    \label{HredDecompD}
    \cbH_r(s)=\cbC_r(s)\cbK_r(s)^{-1}\cbB_r(s)
   + \cbD_r
\end{equation}
where
$\cbC_r(s)\in  \IC^{p \times r}$, 
$\cbB_r(s)\in \IC^{r\times m}$, $\cbD_r\in \IC^{p\times m}$
and 
$\cbK_r(s) \in \IC^{r \times r}$
are obtained by Petrov-Galerkin projection:
%\begin{equation} \label{eq:copred}
%\cbC_r(s) = \cbC(s)\cbV,\quad
%\cbB_r(s) = \cbW^{T}\cbB(s),\quad
%\cbK_r(s) = \cbW^{T} \cbK(s) \cbV,\quad\cbD_r = \cbD.
%\end{equation}
\begin{equation} \label{eq:copred}
\begin{array}{ccc}
\cbK_r(s) = \cbW^{T} \cbK(s) \cbV,  & ~\cbB_r(s) = \cbW^{T}\cbB(s), \\  \\
\cbC_r(s) = \cbC(s)\cbV, & \mbox{and}~~\cbD_r = \cbD.
\end{array}
\end{equation}
\subsection{Interpolatory Projections for Generalized Coprime Factors}
Interpolatory projections for generalized coprime representations of 
transfer functions were introduced in \cite{beattie2008ipm}. We follow the notation in
 \cite{beattie2008ipm} and use ${\mathcal D}_{\sigma}^{\ell}f$  
 to denote the $\ell^{th}$ derivative of the univariate function $f(s)$ 
evaluated at $s=\sigma$, with the usual convention that 
${\mathcal D}_{\sigma}^{0}f=f(\sigma)$.

\begin{thm} \label{GenIntrp}
Given the original model  transfer function $$\cbH(s)=\cbC(s)\cbK(s)^{-1}\cbB(s) 
     + \cbD,$$ let 
let $\cbH_r(s)$ denote the reduced transfer function in (\ref{HredDecompD}) obtained by projection as in (\ref{eq:copred}) using the model reduction bases $\cbV$ and $\cbW$. 
For the interpolation points $\sigma,\,\mu \in \IC$, suppose that $\cbB(s)$, $\cbC(s)$,  and 
   $\cbK(s)$ are analytic at 
   $\sigma \in \IC$ and $\mu \in \IC$. Also let
    $\cbK(\sigma)$ and $\cbK(\mu)$
    have full rank.    Also, let $\rdir\in \IC^{m}$ and  $\ldir\in \IC^{\ell}$  be nontrivial  
    tangential direction vectors. The following implications hold: 
    \begin{itemize}
        \item[(a)]~If 
	\begin{equation}  \label{eq:coprv}
	{\mathcal D}_{\sigma}^{i}[\cbK(s)^{-1}\cbB(s)]\rdir\in 
    \mbox{\normalfont\textsf{Ran}}(\cbV)\quad \mbox{ for } i=0,\ldots,\,N,\end{equation}
    then 
    \begin{equation} 
    \cbH^{(\ell)}(\sigma)\rdir=\cbH_{r}^{(\ell)}(\sigma)\rdir,\quad
     \mbox{for } \ell=0,\ldots,\,N. \end{equation}
        \item[(b)]~If 
	\begin{equation}\label{eq:coplv}
	 \left(\ldir^T{\mathcal D}_{\mu}^{j}
	[\cbC(s)\cbK(s)^{-1}]\right)^{T}\in\mbox{\normalfont\textsf{Ran}}(\cbW)
 \quad \mbox{ for } j=0,\ldots,\,M,\end{equation}
	then 
	\begin{equation}
	\ldir^T\cbH^{(\ell)}(\mu)=\ldir^T\cbH_{r}^{(\ell)}(\mu)\quad
     \mbox{for }\ell=0,\ldots,\,M.\end{equation} \\\vspace{-3ex}
     \item[(c)]~If both $($\ref{eq:coprv}$)$ and $($\ref{eq:coplv}$)$  hold and if $\sigma = \mu$, then 
     \begin{equation}
     \ldir^T\cbH^{(\ell)}(\sigma)\rdir=\ldir^T\cbH_{r}^{(\ell)}(\sigma)\rdir,\quad
      \mbox{for }\ell=0,\ldots,\,M+N-1.
     \end{equation}
    \end{itemize}
        assuming $\cbK_{r}(\sigma)=\cbW^{T}\cbK(\sigma)\cbV$, and
    $\cbK_{r}(\mu)=\cbW^{T}\cbK(\mu)\cbV$ have full rank.
\end{thm}
Theorem \ref{GenIntrp} proves the power and flexibility of the interpolatory framework for model reduction. 
The earlier interpolation result, Theorem \ref{thm:interpolation_highorder},  directly extends to
  this much more  general class of transfer function, requiring very similar subspace conditions. Moreover, 
 this structure guarantees that the reduced transfer function will have a similar generalized coprime representation.
 Computational complexity is comparable; one need only solve $n \times n$ (often sparse) linear systems. 
 
 Recall the delay example (\ref{ExDelaySS}).  Assume that 
 $r=2$ interpolation points $\{\sigma_1,\sigma_2\}$ together with the right-directions
$\{ \rdir_1,\rdir_2\}$  and the left-directions
$\{ \ldir_1,\ldir_2\}$ are given.  Based on Theorem \ref{GenIntrp}, we construct 
$\cbV\in \IC^{n\times 2}$ and $\cbW\in \IC^{n\times 2}$ using
$$ 
\cbV =\left[ \left(\sigma_1\,\bE -\bA_{0} -e^{-\tau \sigma_1}\, \bA_{1}\right)^{-1}\bB \rdir_1\quad\left(\sigma_2\,\bE -\bA_{0} -e^{-\tau \sigma_2}\, \bA_{1}\right)^{-1}\bB \rdir_2 \right] 
$$
and
$$
\cbW=\left[ \left(\sigma_1\,\bE -\bA_{0} -e^{-\tau \sigma_1}\, \bA_{1}\right)^{-T}\bC^T \ldir_1\quad\left(\sigma_2\,\bE -\bA_{0} -e^{-\tau \sigma_2}\, \bA_{1}\right)^{-T}\bC^T \ldir_2 \right].
$$
Notice, first of all, that structure is preserved: The reduced model of dimension-$2$ has the same internal delay structure
\begin{align} \label{ExDelaySSROM}
\cbW^{T}\bE\cbV\,\dot{\bx}_{r}(t)  &= \cbW^{T}\bA_{0}\cbV\,\bx_{r}(t) + 
\cbW^{T}\bA_{1}\cbV\,\bx_{r}(t-\tau) 
+\cbW^{T}\bB\, \bu(t)  \\
 & \quad \by(t) = \bC\cbV\, \bx_{r}(t)\, \nonumber
\end{align}
with a correspondingly structured transfer function
$$
 \cbH_r(s) = \bC\cbV\left(s\,\cbW^{T}\bE\cbV -\cbW^{T}\bA_{0}\cbV -e^{-\tau s}\, \cbW^{T}\bA_{1}\cbV\right)^{-1}\cbW^{T}\bB.
$$
Moreover, due to the interpolation-based construction of $\cbV$ and $\cbW$, the reduced transfer function is a Hermite bitangential interpolant:
$$\cbH(\sigma_i)\rdir_i=\cbH_{r}(\sigma_i)\rdir_i, ~\ldir_i^T\cbH(\sigma_i)=\ldir_i^T\cbH_{r}(\sigma_i)~\mbox{and}~ 
\ldir_i^T\cbH'(\sigma_i) \rdir_i = \ldir_i^T\cbH_{r}'(\sigma_i) \rdir_i$$ for $i=1,2$. Note that the reduced transfer function fully incorporates the delay
structure and {\it exactly} interpolates the original transfer function. This would not be true if the delay term $e^{-\tau s}$ had been approximated via a rational approximation, such as Pad\'{e} approximation, as is commonly done while seeking to convert to a standard first-order form. 
 
 \begin{rem}
The construction of rational interpolants  $\cbH_r(s)=\cbC_r(s)\cbK_r(s)^{-1}\cbB_r(s) 
     + \cbD_r$ with $\cbD_r \neq \cbD$ can be achieved for 
     generalized coprime representations similarly as described in Theorem \ref{GenIntrp2} for standard first-order realizations. 
  For details, we refer the reader to the original source \cite{beattie2008ipm}.
 \end{rem}
 \section{Realization Independent Optimal $\htwo$ Approximation}   \label{sec:tfirka}
   We described \IRKA~in Section \ref{sec:irka} and promoted it as an effective tool for
   constructing at modest cost locally optimal rational $\htwo$ approximations.  
  One may observe, however, that the formulation of \IRKA, as it appears in Algorithm \ref{sucratkry},
   assumes that a  first-order realization for $\bH(s)$ is available: $\bH(s) = \bC(s\bE - \bA)^{-1}\bB$. 
  As we found in Section \ref{sec:coprime}, there are several important examples
   where the original transfer function, $\bH(s)$, is not naturally represented in this way.
In order to address these situations among others, Beattie and Gugercin in \cite{beattie2012realization}
  removed the need for \emph{any} particular realization and extended  applicability of 
  \IRKA~to any evaluable $\htwo$-transfer function. 
   We  focus on this extension of \IRKA~in the present section.  Since $\bH(s)$ is not required to have a first-order realization here, we will follow the notation of Section \ref{sec:coprime} and use $\cbH(s)$  to denote the original  transfer function, however we do not require even that the original system have a generalized coprime representation. Since the reduced model is still a rational function of order $r$ in the first-order form, we continue to use $\bH_r(s)$ to denote the reduced model transfer function.
   
 There are two main observations behind the methodology offered in \cite{beattie2012realization}. The first observation is based on the first-order $\htwo$ optimality conditions  (\ref{H2optcond}) in Theorem \ref{h2cond}. Recall that Theorem \ref{h2cond} does not put any restrictions on $\cbH(s)$; the only assumption is that the approximant  $\bH_r(s)$ is a rational function; thus the theorem and the bitangential Hermite optimality conditions apply equally  if $\cbH(s)$ were to have the form, for example,  
 $\cbH(s) = \bC(s^2\bM + s\bG + \bK)^{-1}\bB$. A second observation is related to how \IRKA~constructs the solution: For  the current set of interpolation points and tangent directions,   \IRKA~constructs a bitangential Hermite interpolant and updates the interpolation data. Thus,  the key issue becomes, given a set of interpolation data, how shall one   
 construct a rational approximant $\bH_r(s)$ that is a Hermite bitangential interpolant   to $\cbH(s)$ (which may not be presented as a first-order state-space model). The Loewner interpolatory framework introduced by Mayo and Antoulas \cite{mayo2007fsg} (discussed in detail in Chapter  \ref{chap:Antoulas}) is the right tool.

 For the balance of this section, we assume that both $\cbH(s)$ and its derivative, $\cbH'(s)$, are only accessible through the evaluation, $s\mapsto (\cbH(s),\cbH'(s))$.  No particular system realizations are assumed. 
 Suppose we are given interpolation points $\{\sigma_1,\ldots,\sigma_r\}$ 
 with the corresponding tangential directions
$\{\rdir_1,\ldots, \rdir_r\}$ and $ \{\ldir_1,\ldots, \ldir_r\}$. We want to construct a degree-$r$ rational approximant $\bH_r(s)$ that is a bitangential Hermite interpolant to $\cbH(s)$:  
\begin{subequations} \label{loewnerconditions}
\begin{align}
 \cbH({\sigma}_k) \rdir_k &=\bH_r({\sigma}_k) \rdir_k,\label{loewnerright} \\
 \rdir_k^T \cbH({\sigma}_k) &=  \rdir_k^T\bH_r({\sigma}_k),  \mbox{ and }
 \label{loewnerleft} \\
 \rdir_k^T \cbH'({\sigma}_k) \rdir_k 
   &=  \rdir_k^T\bH_r'({\sigma}_k) \rdir_k \label{loewnerhermite} 
\end{align}
\end{subequations}
for $k=1,\,2,\,...,\,r$.
As seen in Chapter \ref{chap:Antoulas}, the framework of \cite{mayo2007fsg} allows one to achieve this goal  requiring only the evaluation $\cbH(s)$ and $\cbH'(s)$ at $\sigma_k$ without any constraint on the structure of $\cbH(s)$: Simply construct
the matrices $\bE_r$,  $\bA_r$, $\bB_r$ and $\bC_r$ using
{ 
\begin{equation} \label{defineL}
\left(\bE_r\right)_{i,j} :=
\left\{  \begin{array}{ll}
\displaystyle 
 - \frac{\ldir_i^T\left(\cbH(\sigma_i)-\cbH(\sigma_j)\right)\rdir_j}{\sigma_i-\sigma_j}  & {\rm if}~~ i \neq j \\ \\
- \ldir_i^T  \cbH'(\sigma_i)\rdir_i & {\rm if}~~ i = j 
\end{array}, \right.
\end{equation}
\vspace{2ex}
\begin{equation}\label{defineM}
\left(\bA_r\right)_{i,j} := 
\left\{  \begin{array}{ll}
\displaystyle 
 - \frac{\ldir_i^T\left(\sigma_i\cbH(\sigma_i)-\sigma_j\cbH(\sigma_j)\right)\rdir_j}{\sigma_i-\sigma_j}  &  {\rm if}~~i \neq j \\ \\
- \ldir_i^T  \left. [s\cbH(s)]'\right|_{s=\sigma_i} \rdir_i & {\rm if}~~i= j 
\end{array}, \right.
\end{equation}
}
and 
{\small 
\begin{equation}\label{defineZY}
\bC_r = [\cbH(\sigma_1)\rdir_1,\ldots,\cbH(\sigma_r)\rdir_r],\ 
~\bB_r= \left[ \begin{array}{c} 
\ldir_1^T\cbH(\sigma_1)\\ \vdots \\ \ldir_r^T\cbH(\sigma_r) \end{array}
\right].
\end{equation}
}
Then $\bH_r(s) = \bC_r(s\bE_r-\bA_r)^{-1}\bB_r$ 
  satisfies  (\ref{loewnerconditions}).  
  
To use \IRKA~for $\htwo$ approximation without any structural constraints on $\cbH(s)$, one need only replace the projection-based construction of the intermediate Hermite interpolant with a Loewner-based construction. This is exactly what \cite{beattie2012realization}  introduced, leading to the following realization-independent optimal $\htwo$ approximation methodology:
  \begin{center}
    \framebox[4.5in][t]{
\vspace*{-0.5cm}
    \begin{minipage}[c]{4.3in}
    {\small
    \begin{alg} \label{tfirka} 
    {\bf TF-IRKA:  \IRKA~using transfer function evaluations}
    \begin{enumerate}
       \item Make an initial $r$-fold shift selection: $\{\sigma_1,\ldots,\sigma_r\}$ that is closed under conjugation      
       (i.e., $\{\sigma_1,\ldots,\sigma_r\}\equiv\{\overline{\sigma_1},\ldots, \overline{\sigma_r}\}$ viewed as sets) 
        and initial \\ tangent directions 
       ${\rdir}_1,\ldots,{\rdir}_r$ and 
       ${\ldir}_1,\ldots,{\ldir}_r$, also closed under conjugation.
       \item while (not converged)
       \begin{enumerate}
       \item Construct  $\bE_r$, $\bA_r$, $\bC_r$ and $\bB_r$ as in 
       $($\ref{defineL}$)$-$($\ref{defineZY}$)$.
       \item Compute a pole-residue expansion of $\bH_r(s)$:
         $$
         \bH_r(s) = \bC_r(s \bE_r -\bA_r)^{-1}\bB_r = \sum_{i=1}^r\frac{\ldir_i\rdir_i^T}{s-{\lambda}_i}
         $$
              \item $\sigma_i \longleftarrow -\lambda_i,\quad\rdir_i \longleftarrow \widehat{\rdir}_i,\quad\mbox{ and } 
\quad\ldir_i \longleftarrow \widehat{\ldir}_i,\quad \mbox{for}\quad i=1,\ldots,r$
       \end{enumerate}
           \item Construct  $\bE_r$, $\bA_r$, $\bC_r$ and $\bB_r$ as in 
           $($\ref{defineL}$)$-$($\ref{defineZY}$)$.
         \end{enumerate}
       \end{alg}
       }
    \end{minipage}
    }
    \end{center}
    As for the original formulation of \IRKA, upon convergence the rational approximant resulting from Algorithm \ref{tfirka} will satisfy the first-order necessary conditions (\ref{H2optcond}) for $\htwo$ optimality.

 \subsection{An optimal rational approximation for a delay system}
 Consider the delay system given in (\ref{eq:delayH}), i.e., 
$$
 \cbH(s) = {\bC} {\left(s\,\bE -\bA_{0} -e^{-\tau s}\, \bA_{1}\right)}\,^{-1} {\bB}.
$$
%$$
%\begin{array}{l} 
%\bfE =  \kappa\, \bfI + \bfT \\  \\
%\bfA_0 = \frac{3}{\tau}\left(\bfT-\kappa \,\bfI\right) \\  \\
%\bfA_1 =   \frac{1}{\tau}\left(\bfT-\kappa\,\bfI\right)
%\end{array}\qquad \mbox{with}\quad 
%\bfT=\left[
%\begin{array}{cccccc}
%1  &  1 &  & \ldots  & &  0  \\
%1  &  0 & 1 &           &  &     \\
%    &  1 & 0 &            & &  \\
%\vdots &  &  & \ddots   &   &  \vdots \\
% & & &   &  0 &  1\\
%0 & & &  \ldots  &  1 &  1
%\end{array}
%\right]\in \IR^{n\times n}
%$$
Following \cite{beattie2008ipm}, we take 
$\bE =  \kappa\, \bI + \bT$,
$\bA_0 = \frac{3}{\tau}\left(\bT-\kappa \,\bI\right)$, and
$\bA_1 =   \frac{1}{\tau}\left(\bT-\kappa\,\bI\right)$, for any $\kappa>2$  and delay $\tau>0$,
where $\bT$ is an ${n\times n}$ matrix with ones on the first superdiagonal, on the first subdiagonal, at the $(1,1)$ entry, and at the $(n,n)$ entry. The remaining entries of $\bT$ are zero.
We take the internal delay as $\tau=0.1$ and a  SISO system with $n=1000$, i.e., $\bE,\bA_{0},\bA_1 \in \IR^{1000 \times 1000}$, and $\bB,\bC^T \in \IR^{1000 \times 1}$.  Then, we use \textsf{TF-IRKA} as illustrated in Algorithm \ref{tfirka} to construct a degree $r=20$ (locally) $\htwo$ optimal rational approximation $\bH_r(s)$.  \textsf{TF-IRKA} requires evaluating $\cbH(s)$ and $\cbH'(s)$. For this delay model, $\cbH'(s)$ is given by
$$
 \cbH'(s) = -\bC (s\bE - \bA_1 - e^{-\tau s}\bA_2)^{-1} (\bE + \tau e^{-\tau s}\bA_2) (s\bE - \bA_1 - e^{-\tau s}\bA_2)^{-1}\bB.
$$
Another approach to obtain a rational approximation for such a delay system would be to replace the exponent $e^{-\tau s}$ with a rational approximation and then reduce the resulting rational large-scale model with standard techniques. Here we will use the second order Pad\'{e} approximation towards this goal where we 
replace $e^{-\tau s}$ by 
$\frac{12-6\tau s+\tau^2 s^2}{12+6\tau s+\tau^2 s^2}$, obtaining the large-scale approximate rational transfer function
 \begin{equation} \label{fompade2}
 \cbH^{[P_2]}(s) = (12\bC+s 6\tau\bC+s^2\tau^2\bC) (\bN s^3 + \hat{\bM} s^2 + \hat{\bG} s + \hat{\bK})^{-1}\bB
 \end{equation}
  where 
$
 \bN = \tau^2 \bE$, $
 \hat{\bM} = 6\tau \bE - \tau^2(\bA_0+\bA_1)$,
$\hat{\bG} = 12\bE + 6\tau(-\bA_0+\bA_1)$, and
$\hat{\bK} = -12(\bA_0+\bA_1).$ We used the notation $ \cbH^{[P_2]}(s)$ to denote the resulting large-scale rational approximation due to the second-order  Pad\'{e} approximation. We note that the resulting approximation has a term $s^3$ and will result in an $N=3000$ first-order model. 
Once an interpolatory model reduction technique is applied to $ \cbH^{[P_2]}(s)$, the reduced model will be an exact interpolant to  $ \cbH^{[P_2]}(s)$, but not to $\cbH(s)$. This starkly contrasts to \textsf{TF-IRKA} where the resulting rational approximation exactly interpolates the original delay model 
$\cbH(s)$. In Figure \ref{fig:tfirka_vs_pade} below, we show the Amplitude Bode Plots of $\cbH(s)$
(denoted by ``Full"), order $r=20$ \textsf{TF-IRKA} approximant (denoted by ``TF-IRKA"), and order 
$N=3000$ Pad\'{e} model (denoted by ``Pade"). The figure  illustrates clearly that the locally optimal $\htwo$ approximation due to \textsf{TF-IRKA} almost exactly replicates the original model. On the other hand, even the order $N=3000$ Pad\'{e} model  $\cbH^{[P_2]}(s)$ is a very poor approximation.
\begin{figure}[hh]
   \centerline{
  \begin{tabular}{c}
  \epsfxsize=100mm
   \epsfysize=75mm
      \epsffile{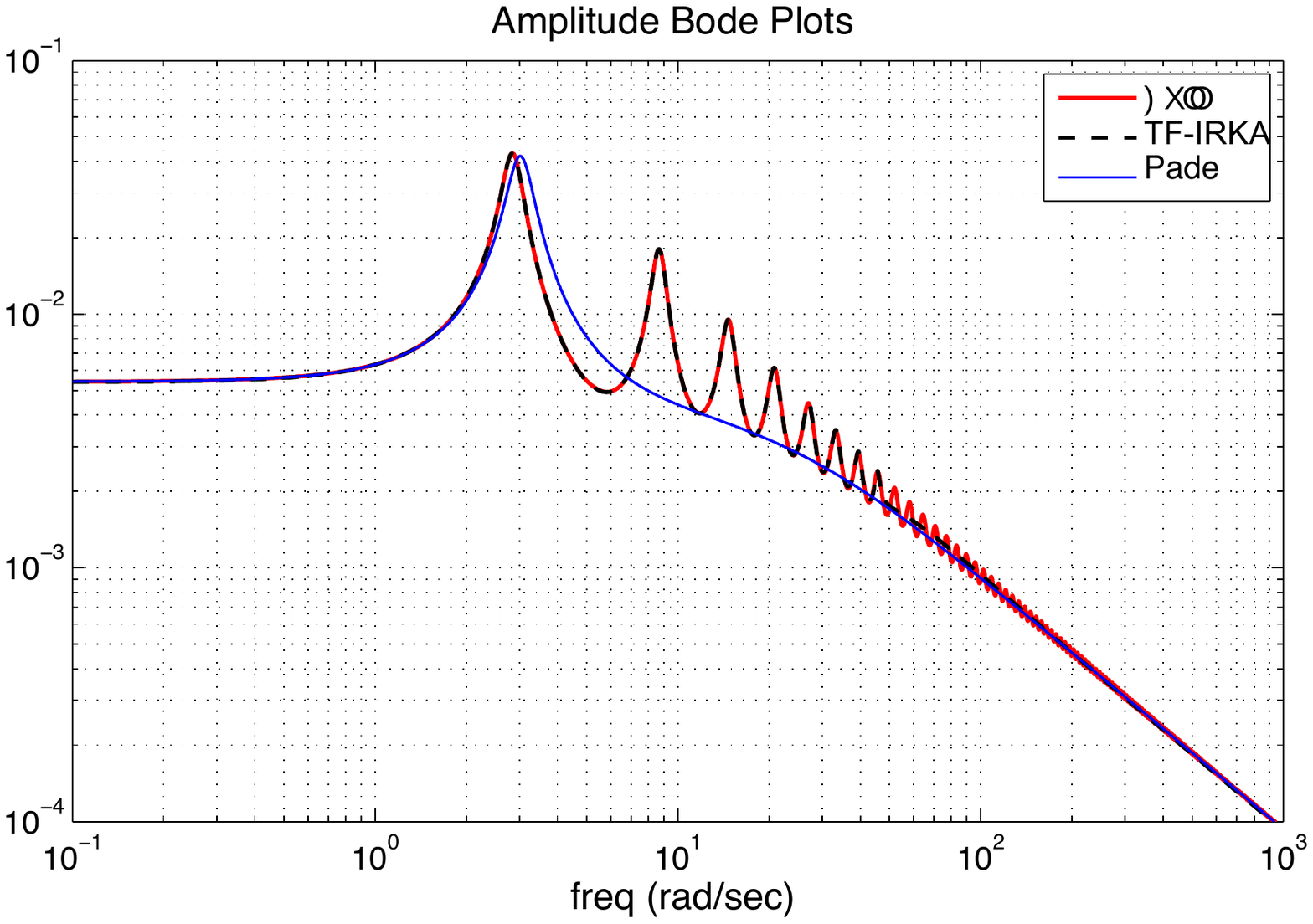}
 \end{tabular}
 } 
 \label{fig:tfirka_vs_pade}
   \end{figure}

\section{Interpolatory Model Reduction of  Parametric Systems}  \label{sec:pmor}
All dynamical systems considered here are \emph{linear, time invariant} systems;  the system properties are presumed to be \emph{constant} at least with respect to time.  Very often system properties will depend on external parameters and system dynamics vary as parameter values change.  Parameters  enter naturally into the models in various ways, representing changes in boundary conditions, material properties,  system geometry, etc.  Producing a new reduced model for every new set of parameter values could be very costly, so a natural 
goal is to generate parametrized reduced models that provide high-fidelity approximations throughout a wide range of parameter values. This is usually referred to as  \emph{parametric model reduction (PMR).} It has found immediate applications in inverse problems \cite{Wang:2005b, Galbally2009, Lieberman2010,Druskin2011solution,morpals2013}, optimization \cite{Arian2002,Antil2011,yue2012,yue2013,Antil2012}, and
design and control  \cite{Lieu2007,
  Amsallem2008,Daniel2004,Baur_etal2011,feng2005pfc,Lieu2006,Hay09}.  There are various approaches to parametric model reduction methods; see, e.g., \cite{prud2002reliable,RozHP08,BuiThanh2008,nguyen2008best,BauB09,haasdonk2011erm,Hay09,veroy2003posteriori,BuiThanh2008_AIAA} and the references therein. In this section, we focus on  interpolatory methods. For a recent, detailed survey on parametric model reduction, we refer the reader to \cite{BGW2013}.
\subsection{Parametric Structure}
We consider  MIMO transfer functions that are parametrized with $\nu$ parameters
 $\parv=[\parv_1,\,\ldots,\, \parv_\nu]$:
 \begin{equation}    \label{paramH}
     \cbH(s, \parv)=\cbC(s, \parv)\cbK(s, \parv)^{-1}\cbB(s, \parv)
\end{equation}
with $\cbK(s, \parv) \in \IC^{n \times n}$ and  $\cbB(s, \parv) \in \IC^{n \times m}$ and  $ \cbC(s, \parv)\in \IC^{p \times n}$. The standard case of parametric linear dynamical systems of the form
\begin{equation} \label{eq:pss}
\bE(\parv)\, \dot{\bx}(t;\parv)  = \bA(\parv)\, \bx(t;\parv) + \bB(\parv)\, \bu(t),
\quad \by(t;\parv)  = \bC(\parv)\, \bx(t;\parv)
\end{equation}
then becomes a special case of the more general form (\ref{paramH}) we consider here with
$\cbK(s, \parv) = s\bE(\parv) - \bA(\parv)$, $\cbB(s, \parv) = \bB(\parv)$ and 
$\cbC(s, \parv)= \bC(\parv)$. 

Even though the theoretical discussion applies to general parametric dependency, we assume an affine parametric form:
\begin{align} 
\cbK(s,\parv)&=\cbK^{[0]}(s)  + k_1(\parv) \, \cbK^{[1]}(s)\, +\ldots+ k_{q}(\parv) \, \cbK^{[\nu]}(s)\,
 \nonumber\\
\cbB(s,\parv)&=\cbB^{[0]}(s)  + b_1(\parv)\, \cbB^{[1]}(s)\ +\ldots+ b_{q}(\parv)
\, \cbB^{[\nu]}(s), \label{paramEABC}\\
\cbC(s,\parv)&=\cbC	^{[0]}(s)  + c_1(\parv)\, \cbC^{[1]}(s)\ +\ldots+ c_{q}(\parv)
\, \cbC^{[\nu]}(s),
\nonumber
\end{align}
\noindent where $\{k_i(\parv)\}$,$\{b_i(\parv)\}$, and $\{c_i(\parv)\}$ for $i=1,\ldots,q$ are scalar-valued nonlinear (or linear) parameter functions.  Even though we have linear dynamics with respect to the state-variable, we allow nonlinear parametric dependency in the state-space representation.  

Our reduction framework remains the same: Use a Petrov-Galerkin projection to construct, in this case, a reduced parametric model. Thus we will pick two model reduction bases  $\cbV \in {\mathbb C}^{n \times r}$ and $\cbW\in {\mathbb C}^{n \times r}$ and obtain the reduced parametric model
\begin{equation}    \label{Hredpar}
    \cbH_r(s,\parv)=\cbC_r(s,\parv)\cbK_r(s,\parv)^{-1}\cbB_r(s,\parv)
\end{equation}
where we use a Petrov-Galerkin projection to obtain the reduced quantities
$\cbC_r(s,\parv) \IC^{p \times r}$, 
$\cbB_r(s,\parv)\in \IC^{r\times m}$,
and 
$\cbK_r(s) \in \IC^{r \times r}$, i.e.,
\begin{equation} \label{eq:copredpar}
\cbC_r(s,\parv) = \cbC(s,\parv)\cbV,\quad
\cbB_r(s,\parv) = \cbW^{T}\cbB(s,\parv),\quad
\cbK_r(s,\parv) = \cbW^{T} \cbK(s,\parv) \cbV.
\end{equation}
Applying (\ref{eq:copredpar}) to the affine parametric structure (\ref{paramEABC}) yields
 \begin{align} 
\cbK_r(s,\parv)&= \cbW^T \cbK^{[0]}(s)  \cbV +
\sum_{i=1}^\nu 
 k_i(\parv) \,\cbW^T \cbK^{[i]}(s) \cbV,\\
\cbB_r(s,\parv)&= \cbW^T\cbB^{[0]}(s)   +
\sum_{i=1}^\nu 
 b_i(\parv)\, \cbW^T\cbB^{[i]}(s), \label{paramEABCred}\\
\cbC_r(s,\parv)&=\cbC^{[0]}(s)  \cbV+\sum_{i=1}^\nu   c_i(\parv)\, \cbC^{[i]}(s) \cbV. 
\nonumber
\end{align}
The advantages are clear.  The affine structure allows fast online evaluation of the reduced model: all the reduced order coefficients matrices can be precomputed and for a new parameter value, only the scalar nonlinear parametric coefficients need to be recomputed. No operation in the original dimension $n$ is required.

\subsection{Interpolatory Projections for Parametric Model Reduction}
The question we want to answer in this section is how to choose $\cbW$ and $\cbV$ so that the reduced parametric model interpolates the original one. The main difference from the earlier cases is that we now have two variables with which to interpolate, namely the frequency variable $s \in \IC$ and the parameter vector $\parv \in \IR^\nu$. Thus, we will require $\cbH_r(s,\parv)$ to (tangentially) interpolate $\cbH(s,\parv)$ at selected $s$ and $\parv$ values. In particular, this will require choosing both frequency interpolation points and parameter interpolation points.

Interpolatory parametric model reduction has been studied in various papers; see, e.g.,~\cite{WeiMGG99,GunKN02, bond2005pmo,   Daniel2004,feng2005pfc,   LeuK05,MosK06,FarHIetal08,FenB08}. These papers focus on matrix interpolation (as opposed to tangential interpolation) and in some cases are restricted to special cases where parametric dependence is allowed only within a subset of state-space matrices. See 
Chapter \ref{chap:BaurBennerHaasdonkHimpeMartiniOhlberger} for a careful of comparison of various parameterized model reduction strategies.  Baur {\it et al.}  \cite{Baur11}  provide quite a general projection-based framework for approaching
structure-preserving parametric model reduction via tangential interpolation. 
Our discussion below follows  \cite{Baur11}  closely.  However, we note that instead of the standard first-order framework (\ref{eq:pss}) that is considered there, we present results for more general \emph{parametrized} generalized coprime representations as in 
(\ref{paramH}). To keep the presentation concise, we only list the zeroth and first order interpolation conditions:

\begin{thm} \label{GenIntrppmor}
Given $\cbH(s,\parv)=\cbC(s,\parv)\cbK(s,\parv)^{-1}\cbB(s,\parv)$, let 
let $\cbH_r(s,\parv)$ denote the reduced transfer function in $($\ref{Hredpar}$)$ obtained by projection as in $($\ref{eq:copredpar}$)$ using the model reduction bases $\cbV$ and $\cbW$. 
For the frequency interpolation points $\sigma,\,\mu \in \IC$ and the parameter interpolation point
$\bpi \in \IR^\nu$, suppose that $\cbB(s,\parv)$, $\cbC(s,\parv)$,  and 
   $\cbK(s,\parv)$ are analytic with respect to $s$ at 
   and $\mu \in \IC$, and are continuously differentiable with respect to $\parv$ in a neighborhood of $\bpi$. Also let
    $\cbK(\sigma,\bpi)$ and $\cbK(\mu,\bpi)$
    have full rank.    Also, let $\rdir\in \IC^{m}$ and  $\ldir\in \IC^{\ell}$  be the nontrivial  tangential directions vectors. Then, 
    \begin{itemize}
        \item[(a)]~If 
	\begin{equation}  \label{eq:vpar}
	\cbK(\sigma,\bpi)^{-1}\cbB(\sigma,\bpi)\rdir\in 
    \mbox{\normalfont\textsf{Ran}}(\cbV)\end{equation}
    then 
    \begin{equation}  \label{eq:pmorrcond}
    \cbH(\sigma,\bpi)\rdir=  \cbH_r(\sigma,\bpi)\rdir. 
    \end{equation}
        \item[(b)]~If 
	\begin{equation}\label{eq:wpar}
	 \left(\ldir^T \cbC(\mu,\bpi)\cbK(\mu,\bpi)^{-1}\right)^{T}\in\mbox{\normalfont\textsf{Ran}}(\cbW),
 \end{equation}
	then 
	\begin{equation}  \label{eq:pmorlcond}
	\ldir^T\cbH(\mu,\bpi)=\ldir^T\cbH_{r}(\mu,\bpi).
     \end{equation} \\ \vspace{-3ex}
     \item[(c)]~If both (\ref{eq:vpar}) and (\ref{eq:wpar})  hold and if $\sigma = \mu$, then 
     \begin{equation}  \label{eq:pmorlrcond}
     \ldir^T\cbH'(\sigma,\bpi) \rdir =\ldir^T\cbH_{r}'(\sigma,\bpi) \rdir
     \end{equation}
     and
   \begin{equation} \label{eq:nablap}
\nabla_{\parv}\ldir^T\cbH(\sigma,\bpi) \rdir =\nabla_{\parv}\ldir^T\cbH_{r}(\sigma,\bpi) \rdir     \end{equation}
    \end{itemize}
    assuming $\cbK_{r}(\sigma,\bpi)=\cbW^{T}\cbK(\sigma,\bpi)\cbV$, and
    $\cbK_{r}(\mu,\bpi)=\cbW^{T}\cbK(\mu,\bpi)\cbV$ have full rank.
\end{thm}
Once again, the basic interpolatory projection theorem extends directly to a more general setting, in this case to  the reduction of parametric systems. Possibly the most important property here is that, as (\ref{eq:nablap}) shows,   interpolatory projection provides matching the parameter sensitivity without ever computing them, i.e. the subspaces $\cbV$ and $\cbW$ do not contain any information about the parameter sensitivity. Nonetheless, the two-sided projection forces a match with this quantity.  Indeed, 
the Hessian with respect to the parameter vector can be matched similarly by adding  more vectors to the subspace;
see \cite{Baur11} for details. 

\subsubsection*{A simple example} 

Consider a mass-spring-damper system where two masses $m_1$ and $m_2$ are connected with a spring-dashpot pair with spring constant $k_2$ and the damping constant $\parv_2$.  Further assume that the mass $m_1$ is connected to ground by another spring-dashpot pair with spring constant $k_1$ and the damping constant $\parv_1$. Also, suppose that a point external force $u(t)$ is applied to $m_1$ and we are interested in the  displacement of the  mass $m_2$.
Let the state-vector $\bx(t) = [\bx_1(t)~\bx_2(t)]$ consists of the displacements of both masses. Then, the corresponding differential equation is given by 
$$
\bM \ddot{\bx} + \bG \dot{\bx} + \bK \bx  = \bb u(t)~~,y(t) = \bc\,\bx(t)
$$
where  $\bb= [1~\,\,~0]^T$, $\bc = [0~\,\,~1]$,
$$
\bM = 
\left[
\begin{array}{cc}
m_1 & 0 \\ 0 & m_2 
\end{array}  \right],
\bG = 
\left[
\begin{array}{cc}
\parv_1 +\parv_2& -\parv_2 \\ -\parv_2 & \parv_2 
\end{array}  \right], \mbox{ and }
\bK = 
\left[
\begin{array}{cc}
k_1 +k_2& -k_2 \\ -k_2 & k_2 
\end{array}  \right].
$$
Let $m_1=m_2=1$, $k_1 = 2$, $k_2 = 2$. Also let the damping constants be parametric and vary as $\parv_1 \in [0.15,0.25]$ and $\parv_2 \in [0.25,0.35]$. Define the parameter vector $\parv = [\parv_1~\parv_2]^T$.
Then, the damping matrix can be written as 
$$
\bG(\parv) = \parv_1  \left[
\begin{array}{cc}
1 & 0 \\ 0 & 0
\end{array}  \right] + \parv_2 \left[
\begin{array}{rr}
1 & -1 \\ -1 & 1
\end{array}  \right] = \parv_1 \cbK_1 + \parv_2 \cbK_2.
$$
Then, the underlying system becomes a parametric dynamical system with a transfer function of the form (\ref{paramH}), i.e. $\cbH(s, \parv)=\cbC(s, \parv)\cbK(s, \parv)^{-1}\cbB(s, \parv)$ with 
\begin{align}
\cbK(s,\parv) &= \underbrace{s^2 \bM + \bK}_{\cbK^{[0]}(s)}+ \underbrace{\parv_1}_{k_1(\parv)}
\underbrace{s\bK_1}_{\cbK^{[1]}(s)} +\underbrace{\parv_2}_{k_2(\parv)}
\underbrace{s\bK_2}_{\cbK^{[2]}(s)} \\  
\cbB(s, \parv) &= \bb = \cbB^{[0]}(s),\quad \mbox{ and }\,\cbC(s, \parv) = \bc =  \cbC^{[0]}(s).
\end{align}
We would like to construct a degree-$1$ parametric reduced model using the frequency interpolation point 
$\sigma= 1$ and the parameter interpolation vector $\bpi = [0.2~0.3]^T$. Note that since the system is SISO, no direction vectors are needed. Then, 
\begin{eqnarray*}
\cbV = \cbK(1,\bpi)^{-1} \cbB(1,\bpi) &=& \left[\begin{array}{r} 2.5661 \times 10^{-1} \\ 1.7885 \times 10^{-1} \end{array}\right],\\
\cbW = \cbK(1,\bpi)^{-T} \cbC(1,\bpi)^T &=& \left[\begin{array}{r} 1.7885 \times 10^{-1} \\ 4.2768 \times 10^{-1} \end{array}\right].
\end{eqnarray*}
This leads to a reduced parametric model $\cbH_r(s, \parv)=\cbC_r(s, \parv)\cbK_r(s, \parv)^{-1}\cbB_r(s, \parv)$  with
\begin{align*}
\cbK_r(s,\parv) &= \left({s^2 \cbW^T \bM \cbV+ \cbW^T\bK\cbV}\right)+ {\parv_1}
{\left(s\cbW^T\bK_1\cbV\right)} +\parv_2 \left(
{s\cbW^T\bK_2\cbV}\right) \\  
\cbB_r(s, \parv) &= \cbW^T\bb,\quad \mbox{ and }\,\cbC(s, \parv) = \bc \cbV.
\end{align*}
One can directly check that  at $\sigma = 1$ and $\bpi = [0.2~0.3]^T$,
$$
\cbH(\sigma, \bpi) = \cbH_r(\sigma, \bpi) =  1.7885 \times 10^{-1},
$$
thus  (\ref{eq:pmorrcond}) holds. Since the system is SISO (\ref{eq:pmorrcond}) and (\ref{eq:pmorlcond})  are equivalent.
Note that 
$$
\cbH'(s, \parv) = -\bc \cbK(s,\parv)^{-1} (2 s \bM + \parv_1 \bK_1 +  \parv_2 \bK_2)  \cbK(s,\parv)^{-1} \bb
$$
and similarly for $\cbH_r'(s, \parv)$. Then, by substituting  $s=\sigma = 1$ and $\parv=\bpi = [0.2~0.3]^T$, we obtain
$$
\cbH'(\sigma, \bpi) = \cbH_r'(\sigma, \bpi) = -2.4814 \times 10^{-1},
$$
thus (\ref{eq:pmorlrcond}) holds. We are left with the parametric sensitivity matching condition (\ref{eq:nablap}). One can directly compute the parameter gradients as 
\begin{align*}
\nabla_{\parv}\cbH(s,\parv) = 
\left[ \begin{array}{c} 
-\bc \cbK(s,\parv)^{-1} (s\bK_1)  \cbK(s,\parv)^{-1} \bb  \\
-\bc \cbK(s,\parv)^{-1} (s\bK_2)  \cbK(s,\parv)^{-1} \bb 
\end{array} \right]
\end{align*}
A direct computation yields that at $s=\sigma = 1$ and $\parv = \bpi = [0.2~0.3]^T$,
$$
\nabla_{\parv}\cbH(\sigma,\bpi) = 
 \nabla_{\parv}\cbH_r(\sigma,\bpi) = 
 \left[ \begin{array}{c} 
 -4.5894 \times 10^{-2} \\
 \phantom{-}1.9349\times 10^{-2}
 \end{array} \right]
$$
As this simple example illustrates, by adding one vector to each subspace, in addition to matching the transfer function and its $s$-derivate, we were able to  match the parameter gradients for free; once again we emphasize that no parameter gradient information was added to the subspaces. However, we still match them by employing a two-sided Petrov-Galerkin projection.

Theorem \ref{GenIntrppmor} reveals how to proceed in the case of multiple frequency and parameter interpolation points.
If one is given two sets of frequency points, $\{ \sigma_i\}_{i=1}^K \in \IC$  and $\{ \mu_i\}_{i=1}^K \in \IC$, the parameter points $\{ \bpi^{(j)}\}_{j=1}^L \in \IC^q$ together with the
right  directions $\{ \rdir_{ij}\}_{i=1,j=1}^{K,L} \in \IC^m$ and the
left directions $\{ \ldir_{ij}\}_{i=1,j=1}^{K,L} \in \IC^p$, compute 
$$
\bv_{ij} =  \cbK(\sigma_i,\bpi^{(j)})^{-1}\cbB(\sigma_i,\bpi^{(j)})\rdir_{i,j}~~{\rm and}~~\bw_{ij} =  \cbK(\mu_i,\bpi^{(j)})^{-T}\cbC(\mu_i,\bpi^{(j)})^T\ldir_{i,j}
$$
for $i=1,\ldots,K$ and $j=1,\ldots,L$ and construct 
\begin{eqnarray*} \nonumber
\cbV&=& [\bv_{11},\ldots,\bv_{1L},\bv_{21},\ldots,\bv_{2L},\ldots,\bv_{K1},\ldots,\bv_{KL}] \in \IC^{n \times (KL)} \\
{\rm and} \hspace{1cm} & &\\
\cbW &=& [\bw_{11},\ldots,\bw_{1L},\bw_{21},\ldots,\bw_{2L},\ldots,\bw_{K1},\ldots,\bw_{KL}] \in \IC^{n \times (KL)} . \nonumber
\end{eqnarray*}
and apply projection as in (\ref{eq:copredpar}).
In  practice, $\cbV$ and $\cbW$ might have linearly dependent columns. In these cases applying a rank-revealing QR or an SVD to remove these linearly independent columns will be necessary and will also help decrease the reduced model dimension.

\begin{rem}
We  have focussed here on a global basis approach to interpolatory parametric model reduction in the sense that we assume that the reduction bases $\cbV$ and $\cbW$ are constant with respect to parameter variation and rich enough to carry global information for the entire parameter space. As for other parametric model reduction approaches, interpolatory model reduction can also be 
formulated with $\parv$-dependent model reduction bases, $\cbV(\parv)$ and $\cbW(\parv)$.  These parameter dependent bases can be constructed in several ways, say by interpolating local bases that correspond to parameter samples $\bpi^{(i)}$.  Such considerations are not specific to interpolatory approaches and occur in other parametric model reduction approaches where the bases might be computed via POD, Balanced Truncation, etc.   Similar questions arise in how best to choose the parameter samples $\bpi^{(i)}$. This also is a general consideration for all parametric model reduction methods. The common approaches such as greedy sampling can be applied here as well. For a detailed discussion of these general issues related to parametric model reduction, we refer the reader to \cite{BGW2013}. We mention in passing that
\cite{Baur11}  introduced  an {\it optimal} joint parameter and frequency interpolation point selection strategy for a special case of parametric systems.
 \end{rem}
 
 \section{Conclusions}
We have provided here a brief survey of interpolatory  methods for model reduction of large-scale dynamical systems. In addition to a 
detailed discussion of basic principles for generic first-order realizations, we have presented an interpolation framework for more general system classes that include generalized coprime realizations and parameterized systems. Reduction of systems of differential  algebraic equations are also discussed. An overview of optimal interpolation methods in the $\htwo$ norm including the weighted case, has also been provided. 

            \bibliographystyle{plain}

      \bibliography{BG_Luminybiblio}

% ------------------------------------------------------------------------
\end{document}